\documentclass[10pt]{article}
\usepackage{latexsym,amsmath,amscd,amssymb,graphics}
\makeatletter
\@addtoreset{figure}{section}
\def\thefigure{\thesection.\@arabic\c@figure}
\def\fps@figure{h,t}
\@addtoreset{table}{bsection}

\def\thetable{\thesection.\@arabic\c@table}
\def\fps@table{h, t}
\@addtoreset{equation}{section}

\makeatother

\begin{document}

\newtheorem{theorem}{Theorem}[section]
\newtheorem{definition}[theorem]{Definition}
\newtheorem{lemma}[theorem]{Lemma}
\newtheorem{remark}[theorem]{Remark}
\newtheorem{proposition}[theorem]{Proposition}
\newtheorem{corollary}[theorem]{Corollary}
\newtheorem{example}[theorem]{Example}
\newtheorem{examples}[theorem]{Examples}

\newcommand{\bfi}{\bfseries\itshape}

\newcommand{\mmu}{\ensuremath{\mathbf{J}^{-1}(\mu)/G_{\mu}}}
\newcommand{\mmuh}{\ensuremath{(\mathbf{J}^{-1}(\mu)\cap
M_{(H)}^{G_{\mu}})/G_{\mu}}}
\newcommand{\nmmuh}{\ensuremath{\mathbf{J}^{-1}(\mu)\cap M_{(H)}^{G_{\mu}}}}
\newcommand{\momu}{\ensuremath{\mathbf{J}^{-1}(\mathcal{O}_{\mu})/G}}
\newcommand{\momuh}{\ensuremath{\mathbf{J}^{-1}(\mathcal{O}_{\mu})\cap 
M_{(H)}/G}}
\newcommand{\nmomuh}{\ensuremath{\mathbf{J}^{-1}(\mathcal{O}_{\mu})\cap 

M_{(H)}}}
\newcommand{\omu}{\ensuremath{\mathcal{O}_{\mu}}}
\newcommand{\ozero}{\ensuremath{\widehat{\mathcal{O}}_{0}}}
\newcommand{\obeta}{\ensuremath{\mathcal{O}}_{\beta}}
\newcommand{\ddto}{\ensuremath{\left.\frac{d}{dt}\right|_{t=0}}}
\newcommand{\ddt}{\ensuremath{\frac{d}{dt}}}
\newcommand{\ddso}{\ensuremath{\left.\frac{d}{ds}\right|_{s=0}}}
\newcommand{\sele}{\ensuremath{_{L}}}
\newcommand{\selehd}{\ensuremath{_{L_H}}}
\newcommand{\selehu}{\ensuremath{_{L^H}}}
\newcommand{\selehdk}{\ensuremath{_{L_K}}}
\newcommand{\selehuk}{\ensuremath{_{L^K}}}
\newcommand{\elek}{\ensuremath{N(K)/K}}
\newcommand{\suhk}{\ensuremath{^{K}}}
\newcommand{\shk}{\ensuremath{_{K}}}
\newcommand{\snkd}{\ensuremath{_{N_K}}}
\newcommand{\snku}{\ensuremath{_{N^K}}}
\newcommand{\sm}{\ensuremath{_{\mu}}}
\newcommand{\seta}{\ensuremath{_{\eta}}}
\newcommand{\sbeta}{\ensuremath{_{\beta}}}
\newcommand{\ele}{\ensuremath{N(H)/H}}
\newcommand{\ene}{\ensuremath{\mathfrak{n}}}
\newcommand{\enesmu}{\ensuremath{\mathfrak{n}_{\mu}}}
\newcommand{\gele}{\ensuremath{\mathfrak{l}}}
\newcommand{\geles}{\ensuremath{\mathfrak{l}^*}}
\newcommand{\elesmu}{\ensuremath{\mathfrak{l}_{\mu}}}
\newcommand{\snu}{\ensuremath{_{\nu}}}
\newcommand{\seme}{\ensuremath{_{m}}}
\newcommand{\smu}{\ensuremath{_{\mu}^{(H)}}}
\newcommand{\smus}{\ensuremath{_{\mu}^{(H)\,*}}}
\newcommand{\som}{\ensuremath{_{\mathcal{O}_{\mu}}}}
\newcommand{\somu}{\ensuremath{_{\mathcal{O}_{\mu}}^{(H)}}}
\newcommand{\somus}{\ensuremath{_{\mathcal{O}_{\mu}}^{(H)\,*}}}
\newcommand{\hmh}{\ensuremath{h_{|_{M_{H}}}}}
\newcommand{\J}{\ensuremath{\mathbf{J}}}
\newcommand{\EJ}{\ensuremath{\mathbf{EJ}}}
\newcommand{\K}{\ensuremath{\mathbf{K}}}
\newcommand{\gk}{\ensuremath{\mathfrak{k}}}
\newcommand{\g}{\ensuremath{\mathfrak{g}}}
\newcommand{\gt}{\ensuremath{\mathfrak{t}}}
\newcommand{\gts}{\ensuremath{\mathfrak{t}^{*}}}
\newcommand{\sgt}{\ensuremath{_{\mathfrak{t}}}}
\newcommand{\gmu}{\ensuremath{\mathfrak{g}_{\mu}}}
\newcommand{\h}{\ensuremath{\mathfrak{h}}}
\newcommand{\m}{\ensuremath{\mathfrak{m}}}
\newcommand{\s}{\ensuremath{\mathfrak{s}}}
\newcommand{\pe}{\ensuremath{\mathfrak{p}}}
\newcommand{\q}{\ensuremath{\mathfrak{q}}}
\newcommand{\qs}{\ensuremath{\mathfrak{q}^{*}}}
\newcommand{\subh}{\ensuremath{_{\mathfrak{h}}}}
\newcommand{\subm}{\ensuremath{_{\mathfrak{m}}}}
\newcommand{\subs}{\ensuremath{_{\mathfrak{s}}}}
\newcommand{\gs}{\ensuremath{\mathfrak{g}^{*}}}
\newcommand{\gmus}{\ensuremath{\mathfrak{g}_{\mu}^{*}}}
\newcommand{\hs}{\ensuremath{\mathfrak{h}^{*}}}
\newcommand{\ms}{\ensuremath{\mathfrak{m}^{*}}}
\newcommand{\subms}{\ensuremath{_{\mathfrak{m}^{*}}}}
\newcommand{\subv}{\ensuremath{_{V}}}
\newcommand{\bd}{\ensuremath{\mathbf{d}}}
\newcommand{\inv}{\ensuremath{^{-1}}}
\newcommand{\suh}{\ensuremath{^{H}}}
\newcommand{\ssuh}{\ensuremath{^{(H)}}}
\newcommand{\sus}{\ensuremath{^{*}}}
\newcommand{\suxi}{\ensuremath{^{\xi}}}
\newcommand{\sucero}{\ensuremath{^{\circ}}}
\newcommand{\scero}{\ensuremath{_{\circ}}}
\newcommand{\szero}{\ensuremath{_{0}}}
\newcommand{\slo}{\ensuremath{_{\lambda_{\circ}}}}
\newcommand{\lo}{\ensuremath{\lambda_{\circ}}}
\newcommand{\schio}{\ensuremath{_{\chi_{\circ}}}}
\newcommand{\chio}{\ensuremath{\chi_{\circ}}}
\newcommand{\sko}{\ensuremath{_{\kappa_{\circ}}}}
\newcommand{\ko}{\ensuremath{\kappa_{\circ}}}
\newcommand{\sh}{\ensuremath{_{H}}}
\newcommand{\ssh}{\ensuremath{_{(H)}}}
\newcommand{\suinf}{\ensuremath{^{\infty}}}
\newcommand{\ngmu}{\ensuremath{N_{G_{\mu}}(H)}}
\newcommand{\ngmuh}{\ensuremath{N_{G_{\mu}}(H)/H}}
\newcommand{\ad}{\ensuremath{{\rm ad}}}
\newcommand{\Ad}{\ensuremath{{\rm Ad}}}
\newcommand{\no}{\ensuremath{\nu_{\circ}}}
\newcommand{\sno}{\ensuremath{_{\nu_{\circ}}}}
\newcommand{\vm}{\ensuremath{V_m}}
\newcommand{\svm}{\ensuremath{_{V_m}}}
\newcommand{\wm}{\ensuremath{W_m}}
\newcommand{\swm}{\ensuremath{_{W_m}}}
\newcommand{\pp}{\ensuremath{\mathbb{P}}}
\newcommand{\rr}{\ensuremath{\mathbb{R}}}

\newsavebox{\savepar}
\newenvironment{boxit}{\begin{lrbox}{\savepar}
\begin{minipage}[b]{15.8cm}}{\end{minipage}\end{lrbox}\fbox{\usebox{\savepar}}}

\makeatletter
\title{{\bf  Relative normal modes for nonlinear
Hamiltonian systems}
\footnote{This paper is dedicated to the memory of Joaqu\'{\i}n
Ortega Lacasa.}}
\author{Juan-Pablo Ortega 
\footnote{Institut Nonlin\'eaire de Nice, UMR 129, CNRS-UNSA,
1361, route des Lucioles, 06560 Valbonne, France.
{\texttt Juan-Pablo.Ortega@inln.cnrs.fr}. }}

\date{}
\makeatother
\maketitle

\addcontentsline{toc}{section}{Abstract}
\begin{abstract}
An estimate on the number of distinct relative periodic orbits 
around a stable relative equilibrium in a Hamiltonian system 
with continuous symmetry is given. This result constitutes a generalization to the Hamiltonian symmetric framework of a classical result by Weinstein and Moser on the existence of periodic orbits in the energy levels surrounding  a stable equilibrium.The estimate obtained is very precise in the sense that it
provides a lower bound for the number of relative periodic orbits 
at each prescribed energy and momentum values neighboring the
stable relative equilibrium in question and with any prefixed
(spatiotemporal) isotropy subgroup. Moreover, it is easily
computable in particular examples.  It is interesting to see how
in our result the existence of non trivial relative periodic
orbits requires (generic) conditions on the higher order terms of
the Taylor expansion of the Hamiltonian function, in contrast with
the purely quadratic requirements of the Weinstein--Moser Theorem,
which emphasizes the highly non linear character of the relatively
periodic dynamical objects. 
\end{abstract}

\section{Introduction}

The search for periodic orbits around non 
hyperbolic equilibria of a Hamiltonian system 
has traditionally been one of the main topics 
in classical mechanics. The best known results 
in this direction are due to 
Liapounov~\cite{liapounoff} and Horn~\cite{horn}, 
who solved the non resonant case. The general 
case was solved only in 1973 by A. Weinstein 
who proved the following theorem~\cite{weinstein 73}:

\begin{theorem}[Weinstein]
\label{weinstein moser theorem}
Let $(M,\,\omega,\,h)$ be a Hamiltonian system 
and let $m\in M$ be an equilibrium of the associated 
Hamiltonian vector field $X_h$ such that $h(m)=0$ and the quadratic form $\bd^2 h(m)$ is definite. Then, 
for each sufficiently small positive $\epsilon$, there 
are at least $\frac{1}{2}\dim M$ geometrically 
distinct periodic orbits of energy $\epsilon$.
\end{theorem}
 
Further extensions of this result due to J.  
Moser~\cite{moser 76} justify why this theorem 
is usually referred to as the {\bfi Weinstein--Moser 
Theorem}. Bartsch~\cite{bartsch 97} has studied 
periodic orbits on the zero level set of the 
Hamiltonian in Moser's generalized version 
of the theorem. In the last two references, the definiteness of the
second variation of the Hamiltonian is not required in the
total space but in a smaller subspace called {\bfi  resonance
space} that will be defined later on in the paper.

In this paper we will be interested in Hamiltonian 
systems endowed with a continuous symmetry. 
More specifically, Hamiltonian systems of the form  $(M,\,\omega,\,G,\,\J :M
\rightarrow \g\sus,\, h:
M \rightarrow \mathbb R)$, where $G$ is a Lie group, with Lie algebra $\g$,
acting properly and canonically on the smooth 
symplectic manifold $(M,\,\omega)$, that encodes the symmetries of the system. We will 
assume that the $G$--action admits an equivariant 
momentum map
$\J:M\rightarrow\g\sus$, where $\g\sus$ denotes 
the dual space of $\g$, and that the Hamiltonian 
function $h$ is $G$--invariant (check for instance with~\cite{fom} for an
introduction to these notions). The generalization of Liapounov's
Theorem to the resonant case carried out in the Weinstein--Moser
Theorem is of great relevance in this setup since
the invariance properties associated to symmetries induce resonances
in many occasions.

The Weinstein--Moser Theorem was adapted 
to this category by Montaldi \textit{et al}~\cite{mrs} 
and later by Bartsch~\cite{bartsch 94}, who obtained sharper estimates. Even though these authors worked in  the symmetric 
framework, their papers still dealt with the search 
of periodic orbits near elliptic equilibria. However,  in the
presence of a continuous symmetry, the  critical elements that
generalize equilibria and  periodic orbits to this category are the
so--called  {\bfi relative equilibria (RE)} and {\bfi relative 
periodic orbits (RPOs)}. Recall that a
relative equilibrium of the $G$--invariant 
Hamiltonian $h$
is a point $m\in M$ such that the integral 
curve $m(t)$ of the
Hamiltonian vector field $X_h$ starting at 
$m$ equals
$\exp(t\xi) \cdot m$ for some $\xi \in \g$, where 
$\exp:\g \rightarrow G$ is the exponential map; 
any such $\xi$ is
called a {\bfi velocity\/} of the relative 
equilibrium. Note that if $m$ has a non--trivial 
isotropy subgroup,
$\xi$ is not uniquely determined. The point 
$m\in M$ is said to be a relative periodic orbit 
of the $G$--invariant Hamiltonian $h$ if there 
is a $\tau>0$ and an element $g\in G$ such that
$F_{t+\tau}(m)=g\cdot F_t(m)$ for any
$t\in\mathbb{R}$,
where $F_t$ is the flow of the Hamiltonian 
vector field $X_h$.  The constant $\tau>0$ is 
called the {\bfi relative period\/} of $m$ and the
group element $g\in G$ its {\bfi phase shift\/}. For historical
reasons, we will occasionally refer to the RPOs that we will find in
this paper as {\bfi  relative normal modes} (we already did so in
the title) given that these solutions generalize to the symmetric
context the
\emph{normal modes} or periodic orbits around equilibria provided
by the Liapounov Center Theorem and by the Weinstein--Moser
estimates.

From the point of view of applications a theorem linking stable relative
equilibria to the existence of RPOs presents certain relevance since stable
relative equilibria are known to appear profusely in most common symmetric
Hamiltonian systems: pendula and oscillators subjected to various 
interactions~\cite{marsden scheurle},
symmetric rigid bodies (free, in the presence  gravity~\cite{lrsm},
or immersed in fluids~\cite{leonard 1997, leonard marsden 1997}),
molecules~\cite{molecules}, point vortices in various phase
spaces~\cite{vortices lim, laurent nonlinearity}, etc.

The search for relative equilibria around stable and unstable relative equilibria 
has been the object of~\cite{singular moser}. The simplest and most
straightforward generalization  of the Weinstein--Moser Theorem to
the symmetric context is  obtained by using {\bfi symplectic
reduction}~\cite{mwr}.  If the point $m\in M$ is such that
$\J(m)=\mu$ is a  regular value of the momentum map $\J$ and the 
coadjoint isotropy subgroup $G_{\mu}$ of $\mu\in\gs$ 
acts freely and properly on the level set $\J\inv(\mu)$, 
then the quotient manifold $\J\inv(\mu)/G_{\mu}$ is a 
symplectic manifold and the dynamics of any $G$--invariant 
Hamiltonian on $M$ drops naturally to Hamiltonian 
dynamics on the reduced manifold $\J\inv(\mu)/G_{\mu}$. 
Moreover, REs and RPOs in $M$ coincide with equilibria 
and periodic orbits in the reduced space, respectively. 
Therefore, if $m$ is a RE such that the Hessian of the 
reduced Hamiltonian at the reduced equilibrium satisfies  
the hypothesis of the Weinstein--Moser Theorem, then 
there are at least $\frac{1}{2}\dim(\J\inv(\mu)/G_{\mu})$ 
geometrically distinct periodic orbits on each energy level in  
this reduced space, that lift to as many geometrically 
distinct RPOs in $M$ with momentum $\mu$. We emphasize that when in the symmetric context we 
talk about {\bfi geometrically distinct objects} we mean that
one cannot be obtained from the other by using the relevant group
action in the problem.

One limitation of this method is that it only allows us to prove the existence of RPOs with the same momentum as the stable relative equilibrium whose existence we use as hypothesis. Additionally, if the regularity assumption on the point $m$ is dropped 
in the previous paragraph, the reduced space 
$\J\inv(\mu)/G_{\mu}$ is not a manifold anymore 
but a  Poisson variety in the sense of~\cite{acg, poisson 
reduction singular}, whose symplectic leaves are the 
singular reduced spaces introduced by Sjamaar, 
Lerman, and Bates in~\cite{sl, bl}. 
See also~\cite{thesis, review}. In principle, the 
procedure described in the previous paragraph can 
still be carried out taking, instead of the entire reduced 
space, the smooth symplectic stratum that contains the 
reduced equilibrium. The main inconvenience of this 
approach is the loss of information  that the restriction 
to the stratum implies. For instance, the stratum could 
reduce to a point, in which case the result would be 
empty of content. However, even when the stratum that
contains the relative equilibrium is not trivial the use of the
Weinstein--Moser  Theorem in it does not give us any information
on the neighboring strata which, as we will see, contain non
trivial relative periodic solutions.  A first step in
overcoming these difficulties has been taken  in~\cite{lerman
tokieda} where  a symplectic version of
the Slice Theorem due to Marle, Guillemin, and Sternberg is used
to establish a relation between the RPOs around a given stable
relative equilibrium and the RPOs around the
corresponding symmetric equilibrium in the slice, always when
certain hypotheses on the coadjoint isotropy of the momentum value
of the relative equilibrium are satisfied. This procedure allows
the application of the    Weinstein--Moser Theorem  on the slice to
produce RPOs of the original system. Nevertheless, this treatment
is not optimal and when the relative equilibrium is just an
equilibrium, these theorems do not provide any information, as far
as RPOs is concerned. In our work we will take an approach inspired
by the so called {\it Smale Program}~\cite{smale}  that consists of
analyzing the orbit spaces resulting from quotienting the level
sets of the conserved quantities in the system by the relevant
group action. This strategy will provide results free from some of
the restrictions in~\cite{lerman tokieda} and containing easy to
compute estimates on the number of RPOs around a given stable
relative equilibrium, classified by their energy, momentum value,
and (spatiotemporal) isotropy subgroups.  The main results are
contained in the following three theorems:
\begin{itemize}
\item Theorem~\ref{theorem 1}: this result gives two different lower bounds for the
number of RPOs with prescribed energy, momentum, and isotropy group in a symmetric
Hamiltonian system, neighboring a stable equilibrium with total isotropy. The
symmetry group is assumed to be compact. The estimates provided are based on
two  different critical point theory tools: Lusternik--Schnirelman category and
Morse theory. The former provides an easy to compute dimensional estimate while
the latter is expressed in terms on an Euler characteristic with respect to
equivariant cohomology that, even though is more difficult to
compute it is, in principle, sharper. This result improves the
study carried out in~\cite{lerman tokieda} since
the main results regarding RPOs in those papers are empty of
content when dealing with an equilibrium with total isotropy. 
\item Theorem~\ref{theorem 1 spatiotemporal}: it provides
estimates on the number of RPOs similar to those in Theorem~\ref{theorem 1} but
this time, the predicted solutions have prescribed isotropy
subgroup not only with respect to the symmetry group $G$ of the
system (referred to as the group  of {\bfi  spatial
symmetries}) but with respect to the group $G \times S ^1 $ (group  of {\bfi 
spatiotemporal symmetries}),
where the circle symmetry comes from putting the system in
normal form. This symmetry, that in principle is not a
feature of the given system, is reflected in the
spatio--temporal symmetry properties of the periods of the solutions
predicted by the theorem. Due to the techniques used in the proofs
and the conclusions obtained Theorem~\ref{theorem 1 spatiotemporal}
is NOT a generalization of Theorem~\ref{theorem 1}: in the proof of
Theorem~\ref{theorem 1} intervenes a transversality argument that
guarantees that all the solutions obtained are genuine RPOs and not
just relative equilibria (that could be considered as trivial
RPOs). This conclusion cannot be drawn from Theorem~\ref{theorem 1
spatiotemporal} given that the subgroups of $G \times S ^1 $
intertwine the $G$  and $S ^1 $--actions via the {\bfi  temporal
character} (this terminology will be introduced later on)
preventing us from making the distinction between RPOs and relative
equilibria (see Remark~\ref{difference between theorems}).
\item Theorem~\ref{rpos around stable re}: it generalizes  Theorem~\ref{theorem 1
spatiotemporal} providing, under certain hypotheses, estimates on the number of
RPOs around a stable relative equilibrium.
\end{itemize}

\medskip

The paper is organized as follows:
\begin{itemize}
\item In Section~\ref{preliminaries} we introduce some preliminary material with the purpose of fixing the notation and of future reference. The expert can skip this section.
\item In Section~\ref{the main result} we present the main results that provide an estimate on the number of RPOs surrounding a given stable symmetric equilibrium at each prescribed energy and momentum values neighboring the equilibrium, and with any prefixed spatial and spatiotemporal isotropy subgroup.
\item In Section~\ref{relative periodic orbits around stable relative equilibria} we use the main results in the previous section and the so called reconstruction equations in order to generalize them to an estimate on the number of RPOs around a genuine stable RE.
\end{itemize}

\section{Preliminaries}
\label{preliminaries}

Throughout the paper we will work in the 
category of symmetric Hamiltonian spaces 
whose objects are Hamiltonian systems 
with symmetry 
$(M,\,\omega,\,G,\,\J :M\rightarrow \g\sus,\, h:
M \rightarrow \mathbb R)$. 
Here  $(M,\,\omega)$ is a symplectic 
manifold on which the  Lie group $G$, 
with Lie algebra $\g$, acts properly, 
canonically, and, moreover, in a globally 
Hamiltonian fashion, that is, it admits an 
equivariant momentum map 
$\J:M\rightarrow\g\sus$, where 
$\g\sus$ is the dual space of $\g$. 
The Hamiltonian function $h$ is always 
assumed to be $G$--invariant. 

\subsection{Proper actions, fixed 
point sets, slices, and normalizers.} 
\label{Proper actions, fixed 
point sets, slices, and normalizers.}

The proofs of the facts 
stated below can be  found in~\cite{bredon, kawakubo,  pa, dk}. The
isotropy subgroups associated to  a proper action are always
compact. Let 
$K$  be  a closed subgroup of $G$. The connected 
components of the sets
\begin{eqnarray*}
M^{K}  &=&\{ z\in M\mid K\subseteq G_{z}\}\\
M_{K}  &=&\{ z\in M\mid K=G_{z}\}
\end{eqnarray*}
are submanifolds of $M$. $M\shk$ is an
open submanifold of $M\suhk$. 
$M\suhk$ is usually referred to as the 
set of $K$--{\bfi fixed points} in $M$ and $M\shk$ as the submanifold of {\bfi isotropy type} $K$. If $M$ is a symplectic manifold, then
(the connected components of) $M_{K}$ and $M^{K}$ are 
symplectic submanifolds of $M$. Proper actions are important since they guarantee the existence of slices
and tubular models: let $m\in M$ and $G_m$ be the isotropy
subgroup of  the element
$m$. The Slice Theorem guarantees the existence of a
$G$--equivariant isomorphism 
$\varphi:G\times_{G _m} B\longrightarrow  U$, where $U$ is a
$G$--invariant open neighborhood of  the orbit $G \cdot m $ such
that
$\varphi[e, 0]= m
$,  and where $B$ is an open
$G _{m}$--invariant  neighborhood of
$0$ in the vector space $T\seme M/T\seme(G\cdot m)$, on which $G _{m}$
acts linearly by $h\cdot(v+T\seme(G\cdot m)):=T _m \Phi _h\cdot v+T\seme(G\cdot m)$.
The set $S$ defined as $S:=\varphi [e,\,B]$, is a  (smooth) {\bfi
slice} at $m$ for the $G$--action on $M$.

Slices have very interesting properties. 
A feature that will
be of particular interest  interest to us is the possibility of
using the slice to locally coordinatize the $G$--space $M$ around
the orbit $G\cdot m$ by means of a local cross--section of $G/G_m$.
More specifically, a local cross--section $\sigma$ of the
homogeneous space $G/G_m$ is a differentiable map
$\sigma:\mathcal{Z}\rightarrow G$, where $\mathcal{Z}$ is an open
neighborhood of $G_m$ in $G/G_m$ such that $\sigma(G_m)={\rm e}$
and $\sigma(z)\in z$, for $z\in\mathcal{Z}$. The Slice Theorem for
proper actions of Palais~\cite[propositions 2.1.2 and 2.1.4]{pa}
guarantees that the map
$F:\mathcal{Z}\times S\rightarrow M$ defined by
$F(z,s):=\sigma(z)\cdot s$, is a diffeomorphism onto an open subset
of $M$ that contains $m$. In Section~\ref{relative
periodic orbits around stable relative equilibria} we will briefly
review a symplectic version of this result.

We now suppose that the symplectic manifold $M$ in question is a symplectic vector space $(V,\omega)$ that 
constitutes a symplectic representation space of $G$. In this case
the $K$--fixed point
space $V^K$ is a symplectic vector subspace of $V$, for any
subgroup $K \subset G $. Recall that any symplectic representation
is globally Hamiltonian with an equivariant momentum map
$\J:V\rightarrow\gs$ associated given by
\[
\langle\J(v),\xi\rangle=\frac{1}{2}\omega(\xi\cdot v,v),\qquad\text{for any}\qquad v\in V,\,\xi\in\g.
\]
The symbol $\xi\cdot v$ denotes the infinitesimal generator at $v$ associated to $\xi\in\g$, and $\langle\cdot,\cdot\rangle$ the natural pairing between the Lie algebra $\g$ and its dual. Let now
$N(K)=\{n\in G\,|\,nKn\inv=K\}$ be the normalizer of $K$ in $G$. 
The globally Hamiltonian $G$--action on $V$ induces globally Hamiltonian actions of   $L:=\elek$ on $V_K$ and $V^K$. Moreover, the $L$--action on  $V_K$ is free. The momentum maps $\J\selehuk:V^K\rightarrow\geles$ and $\J\selehdk:V_K\rightarrow\geles$ associated to these actions are given by 
\[
\J\selehuk(v)=\Xi\sus(\J(v)),\qquad\J\selehdk(v)=\Xi\sus(\J(v)),
\]
where $\Xi\sus:(\mathfrak{k}^\circ)^K\rightarrow\geles$ is the
natural
$\elek$--equivariant isomorphism (see~\cite{thesis, review} for the
details) between the $K$--fixed points in the annihilator of $\h$
in $\gs$ and the dual of the Lie algebra $\geles$ of $\elek$.

\subsection{The resonance space and normal form reduction} 
\label{The resonance space and normal form reduction}

Let 
$(V,\,\omega)$ be a symplectic vector space. 
It is easy to show that there is a bijection between 
linear Hamiltonian vector fields on $(V,\,\omega)$ 
and quadratic forms on $V$. Indeed, if 
$A:V \rightarrow V$ is an infinitesimally 
symplectic linear map, that is, a linear Hamiltonian 
vector field on $(V,\,\omega)$, its corresponding 
Hamiltonian function is given by
\[
Q_A(v):=\frac{1}{2}\omega(Av,\,v),
\qquad\text{for any $v\in V$.}
\]
Also, since $A$ belongs to the symplectic Lie algebra $\mathfrak{sp}(V)$, it admits a unique {\bfi Jordan--Chevalley decomposition}~\cite{humphreys, vanderbauwhede 95} of the form $A=A_s+A_n$, where $A_s\in \mathfrak{sp}(V)$ is semisimple (complex diagonalizable), $A_n\in \mathfrak{sp}(V)$ is nilpotent, and $[A_s, A_n]=0$.
If the quadratic form $Q_A$ is definite, a theorem of Krein~\cite{krein 1950, 
moser 58} guarantees that the associated linear 
Hamiltonian vector field $A$ is semisimple 
(complex diagonalizable) and that all its 
eigenvalues lie on the imaginary axis. Let 
$i\no$ be one of the eigenvalues of  $A$ and 
$T\sno:=\frac{2\pi}{\no}$. We define the 
{\bfi resonance space} $U\sno$ of $A$ with 
{\bfi primitive period} $T\sno$ as 
\[
U\sno:=\ker(e^{A_sT\sno}-I).
\]
The resonance space $U\sno$ has the 
following properties (see \cite{williamson, 
golubitsky stewart 87, vanderbauwhede 95}):
\begin{description}
\item[(i)] $U\sno$ is equal to the direct sum 
of the real generalized eigenspaces of $A$ 
corresponding to eigenvalues of the form 
$\pm ik\no$, with $k\in\mathbb{N}\sus$.
\item[(ii)] The pair $(U\sno,\,\omega|_{U\sno})$ 
is a symplectic subspace of $(V, \omega)$.
\item[(iii)] The mapping $\theta\in S^1\mapsto 
e^{\frac{\theta}{\no} A_s}|_{U\sno}$ generates a symplectic 
$S^1$ linear action  on $(U\sno,\,\omega|_{U\sno})$ with associated momentum map $\J_{S^1}:U\sno\rightarrow\mathbb{R}$ given by $\J_{S^1}=\frac{1}{\no}Q_{A_s}|_{U\sno}$.
\item[(iv)] If $(V, \omega)$ is a symplectic 
representation space of the Lie group $G$ 
and the Hamiltonian vector field $A$ is 
$G$--equivariant (equivalently, the quadratic 
form $Q_A$ is $G$--invariant), then the 
symplectic resonance subspace 
$(U\sno, \omega|_{U\sno})$ is also 
$G$--invariant (this follows from the uniqueness of the Jordan--Chevalley decomposition of $A$, which implies that if $A$ is $G$--equivariant, so is $A_s$). Moreover, the $S^1$ 
and $G$ actions on $(U\sno, \omega|_{U\sno})$ 
commute, which therefore defines a  symplectic 
linear action of $G\times S^1$ on $U_{\no}$.
\item[(v)]\textbf{The normal form reduction} \cite{van der meer 85, van der meer 90, vanderbauwhede 95} Let $(V,\,\omega,\,h_\lambda)$ be a $\lambda$--parameter family ($\lambda\in \Lambda$, where $\Lambda$ is a Banach space) of
smooth $G$--Hamiltonian systems such that for any
$\lambda\in\Lambda$,
$h_\lambda (0)=0$, $\bd h_\lambda (0)=0$, and the $G$--equivariant
infinitesimally symplectic linear map $A:=DX_{h\slo}(0)$ is non
singular and has $\pm i\no$ as eigenvalues. Let
$(U\sno,\,\omega|_{U\sno})$ be the resonance space of $A$ with
primitive period $T\sno$. Then, there exist smooth mappings
$\psi:U\sno\times\Lambda\rightarrow V$ and 
$C^{k+1}$--mapping
$\widehat{h_\lambda}:U\sno\times\Lambda\rightarrow\mathbb{R}$ such
that  $\psi(0,\lambda)=0$, for all $\lambda\in\Lambda$,  
$D_{U\sno}\psi(0,\lo)=\mathbb{I}_{U\sno}$ ($D_{U\sno}$ denotes the
partial Fr\'echet derivative relative to the variable in $U\sno$), 
and, most importantly~\cite[Theorem V.5.17]{van der meer
90},~\cite[Theorem 3.2]{vanderbauwhede 95}, if we stay close enough
to zero in
$U\sno$ and to $\lo\in \Lambda$, then the  $S^1$--relative
equilibria  of the $G\times S^1$--invariant Hamiltonian
$\widehat{h_\lambda}$ are mapped by $\psi(\cdot,\lambda)$ to the
set of periodic solutions of $(V,\,\omega,\,h_\lambda)$ in a
neighborhood of $0\in V$, with periods close to $T\sno$. Hence, in
our future discussion we will substitute the problem of searching
periodic orbits for $(V,\,\omega,\,h_\lambda)$ by that of searching
the $S^1$--relative equilibria of the $G\times S^1$--invariant
family of Hamiltonian systems
$(U\sno,\,\omega|_{U\sno},\,\widehat{h_\lambda})$, that will be
referred to as the {\bfi equivalent system}. 
Given any $k \in \Bbb N $, the equivalent system Hamiltonian 
$\widehat{h_\lambda}$ can be chosen so that its Taylor expansion
coincides with that of 
$h_\lambda| _{U\sno}$ up to order $k$. This fact
and the
properties of  $\psi$ imply that
\begin{equation}
\label{restriction to resonance}
\mathcal{A}:=A|_{U\sno}=D\subv X_{h\slo}(0)|_{U\sno}=
D_{U\sno} X_{h\slo|_{U\sno}}(0)=D _{U\sno} X_{\widehat{h\slo}}(0).
\end{equation}
The reduction procedure that we just described appears also in the
literature under the name of {\bfi  Weinstein--Moser reduction}.
\end{description}

\subsection{Hamiltonian relative equilibria and relative periodic orbits} 

A point $m\in M$ is a relative equilibrium of the Hamiltonian system with symmetry $(M,\omega,h, G,\J)$, with velocity $\xi\in\g$, iff $m$ is a critical point of the {\bfi augmented Hamiltonian} $h\suxi:=h-\J\suxi$, where $\J\suxi:=\langle\J,\xi\rangle$. A similar characterization for RPOs that also uses the momentum map is given in the following elementary result.

\begin{proposition}
\label{elementary result}
Let $(M,\,\omega,\,h)$ be a Hamiltonian system with a globally 
Hamiltonian symmetry given by the canonical action of the Lie group
$G$ on $M$, with associated momentum map $\J:M\rightarrow\gs$. If
the Hamiltonian vector field $X_{h-\J\suxi}$, $\xi\in\h$, has a
periodic point $m\in M$ with period $\tau$, then the point $m$ is a
RPO of $X_h$ with relative period $\tau$ and phase shift
$\exp\,\tau\xi$.
\end{proposition}

\noindent\textbf{Proof} Let $F_t$ be the flow of the Hamiltonian vector field $X_h$ and $K_t(m)=\exp\,t\xi\cdot m$ that of $X_{\J\suxi}$. 
By Noether's Theorem we have that $
[X_h,X_{\J\suxi}]=-X_{\{h,\J\suxi\}}=0$,
where the bracket $\{\cdot,\cdot\}$ denotes the Poisson bracket
associated to the symplectic form $\omega$. Due to this equality,
we can write (see for instance~\cite[Corollary 4.1.27]{mta}) the
following expression for $G_t$, the flow of $X_{h-\J\suxi}$,
$
G_t(m)=\lim_{n\rightarrow\infty}(F_{t/n}\circ K_{-t/n})^n(m)=
(K_{-t}\circ F_t)(m)=\exp\,-t\xi\cdot F_t(m)$.
Since by hypothesis the point $m$ is periodic for $G_t$ with period
$\tau$, we have that
$
m=\exp\,-\tau\xi\cdot F_\tau(m)$,
or, equivalently,
$
F_{\tau}(m)=\exp\,\tau\xi\cdot m$,
as required.\ \ \ $\blacksquare$

\medskip

The proposition that we just proved allows us to search the RPOs of
the system $(M,\omega,h, G,\J)$ by looking at the periodic orbits
of  the systems with Hamiltonian functions the augmented
Hamiltonians $h\suxi:=h-\J\suxi$. Notice that in terms of symmetry
properties, the new systems whose periodic orbits we want to
compute are weaker. More specifically, even though the original
Hamiltonian is $G$--invariant, the augmented Hamiltonian
$h-\J\suxi$ is only $G\suxi$--invariant, where $G\suxi$ is the
adjoint isotropy of the element $\xi\in\g$, that is, $G\suxi:=\{g
\in G \mid
\operatorname{Ad}_g \xi = \xi\}$.

As we already mentioned in the introduction, the use of the previous proposition in the search for the RPOs of a system carries intrinsically two main limitations. Firstly, since the phase shift of a RPO that amounts to a periodic orbit of $h-\J\suxi$ is always of the form $\exp\tau\xi$, with  $\tau$ some real
number. Hence, the RPOs whose phase shifts do not lie in the
connected component of the identity of the group of symmetries $G$
cannot possibly be found in this way. Second, if $G$ is a discrete
group then its Lie algebra is trivial and consequently so is the
momentum map associated to this action, which makes the previous
proposition empty of content.

\subsection{Results on critical point theory of functions on 
compact manifolds} 
\label{Results on critical point theory of functions on 
compact manifolds}

\subsubsection{The Lusternik--Schnirelman approach}

The following two results are slight generalizations of those presented in~\cite{weinstein 77} for circle actions. The additional hypotheses that we will introduce in our statements will make the original proofs work with straightforward modifications.

\begin{proposition}
\label{v manifolds}
Let $M$ be a compact $G$--manifold, with $G$ a Lie group acting properly on $M$. Any $G$--invariant smooth function $f\in C\suinf(M)^G$ has at least
\begin{equation}
\label{category bound}
{\rm Cat} (M/G)
\end{equation}
critical orbits. 
\end{proposition} 

In the previous statement, the symbol Cat denotes the {\bfi Lusternik--Schnirelman category} of the quotient compact topological space $M/G$ (the action of $G$ on $M$ does not need to be free and, consequently, the quotient $M/G$ is not in general a manifold). Recall that the  Lusternik--Schnirelman category of a 
compact topological space $M$ is the minimal number of closed
contractible sets needed to cover $M$.  

The preceding  statement is one of  the main analytical tools
that we will use to obtain estimates on the number of RPOs of
our problem. As we will see later on, we will be able to reduce
the search for those RPOs to the computation of the number of
critical orbits of a
$G$--invariant function defined on a compact $G$--symmetric
manifold that satisfies the hypotheses of Proposition~\ref{v
manifolds}. 
In principle, the category in~(\ref{category bound})  is very
difficult to compute. Nevertheless, in our situation we will take
advantage of the symplectic nature of our setup to estimate it in terms
of readily computable dimensional quantities. The main tools
to carry that out are the following technical propositions
whose relevance will become apparent at the time of their
application in the proof of the main theorems.

\begin{proposition}
\label{result b v manifolds}
Let $M$ be a compact $G$--manifold, with $G$ a Lie group acting properly on $M$ such that the isotropy subgroup of each point $m\in M$  is a finite subgroup of $G$. Let $\omega$ be a symplectic $G$--invariant  form defined on $M$. 
Let $H\subset G$ be a Lie subgroup of $G$ and $N\subset M$ be a
$H$--invariant closed submanifold of $M$ such that for any $n\in N$
we have that 
\begin{equation}
\label{reduction hypotheses}
(T_n N)^\omega=\g\cdot n\qquad\text{and}\qquad T_n N\cap\g\cdot m=\h\cdot  n,
\end{equation}
where $\g\cdot n=\left\{\xi_M(n):=\ddto\exp t\xi\cdot n\mid\xi\in\g\right\}$ denotes the tangent space at $n\in N\subset M$ of the $G$--orbit $G\cdot n=\{g\cdot n\mid g\in G\}$.
Then, there is a cohomology class $\theta\in H^2(N/H;\mathbb{R})$ such that $\theta^k\neq 0$, where $k=\frac{1}{2}(\dim N-\dim H)$.
\end{proposition}

\noindent The proof of the following elementary fact can be found in~\cite{schwartz 69}.

\begin{proposition}
\label{cup length}
Let $M$ be a compact topological space. The Lusternik--Schnirelman category of $M$ is at least 1 plus its cuplength.
\end{proposition}

\begin{corollary}
\label{cup length symplectic}
Let $(M,\omega)$ be a $2n$--dimensional compact symplectic manifold. Then,
\begin{equation}
\label{symplectic category}
{\rm Cat} (M)\geq 1+n.
\end{equation} 
\end{corollary}

\noindent\textbf{Proof} The symplecticity of $\omega$ implies that $\omega^n$ is a nowhere vanishing multiple of the volume form, hence $[\omega]^n=[\omega^n]\neq 0$ in the top cohomology group of the manifold, and therefore the cuplength of the manifold $M$ is at least $n$. The conclusion follows from Proposition~\ref{cup length}.\ \ \ \ $\blacksquare$

\medskip

Another approach to the search of critical orbits of symmetric functions is the use of the so called {\bfi equivariant  Lusternik--Schnirelman category} or $G$--Lusternik--Schnirelman category (denoted by the symbol $G$--Cat), introduced in different versions and degrees of generality by Fadell~\cite{fadell 1985}, Clapp and Puppe~\cite{clapp puppe 1986, clapp puppe 1991}, and Marzantowicz~\cite{Marzantowicz 1989}. The equivariant  Lusternik--Schnirelman category is not  the standard Lusternik--Schnirelman category of the orbit space that we used in the previous paragraphs, but the minimal cardinality of a covering of the $G$--manifold $M$ by $G$--invariant closed subsets that can be equivariantly deformed to an orbit. This new category is also a lower bound for the number of critical orbits of a $G$--invariant function on $M$ and is actually a better bound since it can be proven (see for instance~\cite[page 43]{fadell 1985}) that 
$
G-{\rm Cat}(M)\geq {\rm Cat} (M/G)$,
where the equality holds, for instance, when the $G$--action on $M$ is free. Nevertheless, we will not use this category since the cohomological estimates that can be made via arguments similar to those established using propositions~\ref{result b v manifolds} and~\ref{cup length} on the value of the $G$--category require the use of $G$--equivariant cohomology hence giving rise to estimates that are not as readily computable as those that we will obtain using standard cohomology.

\subsubsection{The Morse theoretical approach}

The following paragraphs briefly summarize some results on Morse
Theory adapted to the symmetric Hamiltonian setup as they can be
read in the works of Kirwan~\cite{kirwan 84, kirwan vortices}. 

Let $f\in C\suinf(M)$ be a smooth function on the compact manifold $M$. A critical point $m\in  M$ is said to be {\bfi non--degenerate} when the second derivative of $f$ at the critical point $\bd^2 f(m)$ is a non degenerate quadratic form. Non degenerate critical points are isolated. The {\bfi index} $i_f(m)$ of a critical point $m$ is the maximal dimension of a subspace of $T_mM$ in which $\bd^2 f(m)$ is negative definite. A smooth function all whose critical points are non degenerate is said to be of {\bfi Morse type}.

If $M$ is a $G$--manifold with $G$ a Lie group of positive 
dimension, the critical points of any smooth $G$--invariant
function $f\in C\suinf(M)^G$ come in $G$--orbits. Therefore, since
the critical points of $f$ are not isolated the function $f$ cannot
possibly be of Morse type. The closest condition to being of Morse
type that we can envision in the equivariant context is  what we
will call being of  $G$--Morse. A function $f\in C\suinf(M)^G$ is
of {\bfi $G$--Morse} when all its critical points $m\in M$ satisfy
that
\[
\ker \bd^2 f(m)=T_m(G\cdot m).
\]
This condition will be needed in the statement of our main results. Note that this is a "reasonable" condition to be imposed since it is generic in the $G$--category. A situation in which it is very easy to see that this is the case is when the $G$--action on $M$ is free. In that case, a 
$G$--invariant function being $G$--Morse is equivalent to
its projection onto the quotient manifold $M/G$ being a standard
Morse function, which occurs generically.

The $G$--Morse condition implies the {\bfi minimal degeneracy 
condition}~\cite[Appendix 10]{kirwan 84}, which is the weakest
condition that, to the knowledge of the author, allows the
formulation of Morse inequalities and standard Morse theory
(see~\cite{atiyah bott 82, mcduff salamon} for additional
information). More specifically, the Morse
inequalities in this situation state that if the function $f\in
C\suinf(M)^G$ is $G$--Morse then there exits a polynomial $R(t)$ in
$t$ with non negative integer coefficients such that
\begin{equation}
\label{morse bott inequality}
\sum_{G\cdot x\ \text{critical orbit of}\ f} t^{i_f(G\cdot x)}P_t(G\cdot x)-P_t(M)=R(t)(1+t),\qquad\quad R(t)\geq 0,
\end{equation}
where the sum runs over all the critical orbits $G\cdot x$ of the function $f$, and $P_t(G\cdot x)$, $P_t(M)$ denote the Poincar\'e series of $G\cdot x$ and $M$ respectively
\[
P_t(G\cdot x)=\sum_{i\geq 0}b_i(G\cdot x)t^i,\qquad P_t(M)=\sum_{i\geq 0}b_i(M)t^i.
\] 
The symbols $b_i(M)$ (resp. $b_i(G\cdot x)$) denote the Betti numbers of the manifold $M$ (resp. $G\cdot x$), that is,
\[
b_i(M)=\dim H_i(M)\qquad \text{(resp. }b_i(G\cdot x)=\dim (G\cdot x)\ \text{)}.
\]
The Morse inequalities~(\ref{morse bott inequality}) still hold if instead of using ordinary cohomology we use equivariant cohomology, that is, there exits a polynomial $R(t)$ in $t$ with non negative integer coefficients such that
\begin{equation}
\label{morse bott inequality equivariant}
\sum_{G\cdot x\ \text{critical orbit of}\ f} t^{i_f(G\cdot x)}-P_t^G(M)=R(t)(1+t),\qquad\quad R(t)\geq 0,
\end{equation}
where $P_t^G(M)=\sum_{i\geq 0}t^i\dim H_G^i(M)$. These are the so called {\bfi equivariant Morse inequalities}.
 
A straightforward consequence of ~(\ref{morse bott inequality equivariant}) is that any $G$--invariant function $f$ on the compact manifold $M$ has at least 
\begin{equation}
\label{euler estimate}
|\chi(M)^G|:=|P_{-1}^G(M)|
\end{equation}
critical orbits. The number $\chi(M)^G$ is called the {\bfi $G$--Euler characteristic} of $M$.

The previous remark is particularly relevant in the globally Hamiltonian framework: suppose now that the compact $G$--manifold $M$ is symplectic, that the Lie group $G$ acts canonically, and that this action has an equivariant momentum map associated $\J:M\rightarrow\gs$. Let $\mu\in\J(M)\subset \gs$ be a regular value of $\J$ and $G\sm$ be the corresponding coadjoint isotropy subgroup which we will assume acts in a locally free fashion on $\J\inv(\mu)$.  In the proofs of our main results it will be necessary to evaluate the number of critical orbits of a $G\sm$--invariant function on $\J\inv(\mu)$, which by~(\ref{euler estimate}) can be done just by computing $\chi(\J\inv(\mu))^{G\sm}$. This work has been carried out in full detail by 
Kirwan~\cite{kirwan 84, kirwan vortices} who realized that the
function $f\sm:=\|\J-\mu\|^2\in C\suinf(M)^{G\sm}$ ($\|\cdot\|$ is
any $\Ad_G\sus$--invariant inner product  on $\gs$) is an
equivariantly perfect Morse function, which allows to explicitly
write down (see~\cite[Theorem 4.14]{kirwan vortices}) the
$G\sm$--equivariant Betti numbers of $\J\inv(\mu)$, and
consequently $\chi(\J\inv(\mu))^{G\sm}$, in terms of the Betti
numbers of $M$, of the classifying space of $G\sm$, and of the
equivariant Betti numbers of some simpler subspaces of $M$
explicitly constructed using the geometry of the momentum map.
Therefore, in our future discussions we will consider
$\chi(\J\inv(\mu))^{G\sm}$ as a computable quantity that we can
use, via~(\ref{euler estimate}), in our estimations on the number
of $G\sm$--critical orbits on $\J\inv(\mu)$ (see Section~3
in~\cite{kirwan vortices} for an explicit example of this
calculation).

We conclude this remark by noting that since we will be working with locally free actions, the rational equivariant cohomology ring $H\sus_{G\sm}(\J\inv(\mu);\mathbb{Q})$ is isomorphic to the ordinary rational cohomology ring $H\sus(\J\inv(\mu)/G\sm;\mathbb{Q})$, and therefore
\[
\chi(\J\inv(\mu))^{G\sm}=\chi(\J\inv(\mu)/G\sm).
\]

\subsection{Lagrange multipliers}

The use of Lagrange multipliers will be crucial in the proof techniques that we will use. The version of this result that we present in the following proposition can be found as Corollary 3.5.29 in~\cite{mta}.

\begin{proposition}
\label{lagrange multipliers}
Let $M$ be a smooth manifold, $F$ be a Banach space, $g:M\rightarrow F$ be a smooth submersion, $f\in C\suinf(\mathbb{R})$, and $N=g\inv(0)$. The point $n\in N$ is a critical point of the restriction $f|_N$ iff there exists $\lambda\in F\sus$, called a {\bfi Lagrange multiplier}, such that $n$ is a critical point of $f-\lambda\circ g$. 
\end{proposition}

\section{The main result}
\label{the main result}

\subsection{Relative periodic orbits with prescribed isotropy}
Before we state the main result 
of this section we need some terminology. Let $h\in C\suinf(V)$ a
smooth function defined on the  vector space $V$ such that
$h(0)=0$, $\bd h(0)=0$, and the second derivative at zero $Q:=\bd^2
h(0)$ is a definite quadratic form. Let $\langle\cdot,\cdot\rangle$
be the scalar product on $V$ defined by
\[
\langle u,v\rangle:=\frac{1}{2}\bd^2 h(0)(u,v),\qquad u,v\in V,
\]
and $\|\cdot\|$ be the associated norm. We will say that these are the {\bfi scalar product and the norm associated to the quadratic form} $Q$. We now write the Taylor expansion of $h$ around the origin:
\[
h(v)=\|v\|^2+\frac{1}{3}\bd^3 h(0)\left(v^{(3)}\right)+\cdots+\frac{1}{k!}\bd^k h(0)\left(v^{(k)}\right)+\cdots\qquad v\in V.
\]
We will say that the $k$th term in the Taylor expansion of the function $h$ is {\bfi purely radial} when 
\[
\frac{1}{k!}\bd^k h(0)\left(v^{(k)}\right)=c_k\|v\|^k,
\]
where $c_k$ is a constant real number.

\medskip

The main goal of this section is proving the following theorem:

\begin{theorem}
\label{theorem 1}
Let $(V,\omega,h,G,\J:V\rightarrow\gs)$ be a Hamiltonian system with symmetry, with $V$ a vector space, and $G$ a compact positive dimensional Lie group that acts on $V$ in a linear and canonical fashion. Suppose that $h(0)=0$, $\bd h(0)=0$ (that is, the Hamiltonian vector field $X_h$ has an equilibrium at the origin) and that the linear Hamiltonian vector field $A:=DX_h(0)$ is non degenerate and contains $\pm i\no$ in its spectrum. Let $U\sno$ be the resonance space of $A$ with primitive period $T\sno:=\frac{2\pi}{\no}$.
Let $K\subset G$ be an isotropy subgroup of the $G$--action on $V$ for which the quadratic form $Q\suhk$ on the $K$--fixed point space $U\sno\suhk$ defined by
\[
Q\suhk(v):=\frac{1}{2}\bd^2 h(0)(v,v),\qquad v\in U\sno\suhk
\]
is definite, and the set $\J\selehdk\inv(0)\cap Q\shk\inv(1)$ is non empty ($Q\shk:=Q\suhk|_{(U\sno)\shk}$). In these conditions, there exists a neighborhood $B(K)\subset \geles$ of zero in $\geles$ such that for any $\lambda\in B(K)$, the intersection $\J\selehdk\inv(\lambda)\cap Q\shk\inv(1)$ is a submanifold of $(U\sno)\shk$. Suppose that the following  two generic hypotheses hold:
\begin{enumerate}
\item[{\bf (H1)}]  The restriction $h|_{U\sno\suhk}$ of the Hamiltonian $h$ to the fixed point subspace $U\sno\suhk$ is not radial with respect to the norm associated to $Q\suhk$.
\item[{\bf (H2)}] Let  $h_k(v):=\frac{1}{k!}\bd^k h(0)
\left(v^{(k)}\right)$, $v\in U\sno\suhk$ be the first non radial
term in the Taylor expansion of $h|_{U\sno\suhk}$ around zero. We
will assume that $k\geq 4$ and that the restrictions
$h_k|_{\J\selehdk\inv(\lambda)\cap Q\shk\inv(1)}$ of $h_k$ to the
submanifolds $\J\selehdk\inv(\lambda)\cap Q\shk\inv(1)$, with
$\lambda\in B(K)$, are  $(N(K)/K)_\lambda\times
S^1$--Morse.  
\end{enumerate}
Then, the neighborhood $B(K)$  can be chosen so that for any  $\epsilon>0$ close enough to zero and any $\lambda\in B(K)$ there are at least
\begin{equation}
\label{main estimate 1}
{\rm max}\left[\frac{1}{2}\left(\dim U\sno^K-\dim (N(K)/K)-\dim\left(N(K)/K\right)_{\lambda}\right),\chi\left(\J\selehdk\inv(\lambda)\cap Q\shk\inv(1)\right)^{L_\lambda\times S^1}\right]
\end{equation}
distinct relative periodic orbits of $X_h$ with energy $\epsilon$, 
momentum $\epsilon(\Xi\sus_K)\inv(\lambda)\in\gs$, isotropy
subgroup $K$, and relative period close  to $T\sno$. The
choice of $B (K) $ guarantees that these relative periodic orbits
are not just relative equilibria.  The symbol $N(K)$ denotes the
normalizer of
$K$ in
$G$, 
$\J_{L_K}:(U\sno)_K\rightarrow\geles$ is the momentum map
associated to the free $L:=N(K)/K$--action on $(U\sno)_K$,
$\Xi\sus_K$ denotes the natural isomorphism
$\Xi\sus_K:(\h^\circ)^K\rightarrow\geles$, and 
$\left(N(K)/K\right)_{\lambda}$ the coadjoint isotropy of
$\lambda\in\geles$ (see Section~\ref{Proper actions, fixed  point
sets, slices, and normalizers.} for more information about this
notation). The symbol $\chi\left(\J\selehdk\inv(\lambda)\cap
Q\shk\inv(1)\right)^{L_\lambda\times S^1}$ denotes the
$L_\lambda\times S^1$--Euler characteristic of 
$\J\selehdk\inv(\lambda)\cap Q\shk\inv(1)$ (which in this case
equals the standard Euler characteristic of the symplectic 
quotient 
$\chi(\J\selehdk\inv(\lambda)\cap Q\shk\inv(1)/L_\lambda
\times S^1)$).
\end{theorem}

\begin{remark}
\normalfont
The main estimate~(\ref{main estimate 1}) contains two parts. The
first one, in terms of various dimensions, is in general easy to
compute (see the examples below) and has been obtained via a
cohomological estimate on the Lusternik--Schnirelmann category of
certain symplectic orbifold. The second part involving the Euler
characteristic is in general a sharper estimate but also much more
difficult to compute. The reader is referred to~\cite{kirwan
vortices} for an example on how to compute these quantities.\ \ \
$\blacklozenge$
\end{remark}

\begin{remark}
\normalfont
As we will see in the proof, the properties  of the non
trivial
$S^1$--action on the space $U\sno^K$ (normal form) imply that the
number $k$ in the statement of hypothesis~\textbf{(H2)} is
necessarily even.\ \ \ $\blacklozenge$
\end{remark}

\begin{examples}
\label{examples equilibria}
\normalfont
We illustrate with a few elementary examples the use of
Theorem~\ref{theorem 1}  in specific situations. 
\begin{description}
\item [(i)] {\bf Nonlinearly perturbed spherical pendulum:} the spherical pendulum 
consists of a particle of mass $m$, moving under
the action of a constant gravitational field of
acceleration $g$, on the surface of a sphere of
radius $l$. We consider this system perturbed by
an $S ^1$--invariant nonlinear term. If we use as local coordinates of
the configuration space around the downright position the Cartesian
coordinates
$(x,y)$ of the orthogonal the projection of the
sphere on the equatorial plane, the (local) Hamiltonian of
this system is:
\[
h(x,y,p_x,p_y)=\frac{p_x^2}{2m}+
\frac{p_y^2}{2m}-\frac{(xp_x+yp_y)^2}{2ml^2}
-mg\sqrt{l^2-x^2-y^2}+ \varphi(x ^2+ y ^2, p _x^2+ p _y^2, x p _x + y
p _y),
\]
where the function $\varphi$ is of order two or higher in all of its
variables and encodes the nonlinear perturbation. This system is
invariant with respect to the globally Hamiltonian $S^1$--action
given by the expression
$\Phi_\theta(x,y,p_x,p_y)=
(R_\theta(x,y),R_\theta(p_x,p_y))$, where
$R_\theta$ denotes a rotation of angle $\theta$.
The momentum map
$\J:\mathbb{R}^4\rightarrow\mathbb{R}$ associated
to this action is given by
$\J(x,y,p_x,p_y)=xp_y-yp_x$. The point
$(x,y,p_x,p_y)=(0,0,0,0)$ is an stable equilibrium of the
Hamiltonian vector field $X_h$ to which we
apply Theorem~\ref{theorem 1}. Indeed, the linearization of $X _h$
at the origin has two imaginary eigenvalues that are forced to be
double by the symmetry of the problem. Consequently, the associated
resonance space coincides with the entire $\Bbb R^4$. We now use
the dimensional estimate in~(\ref{main estimate 1}) to look for
RPOs with trivial isotropy. A straightforward dimension count shows
that when $K =\{e\} $, the orbifold $\J\inv(0)\cap
Q_{e}\inv(1)/S ^1
\times S^1 $ is zero dimensional and therefore hypotheses {\bf (H1)}
and {\bf (H2)} trivially hold for the standard spherical pendulum
($\varphi\equiv 0$) and for ALL non radial perturbations $\varphi $.
Therefore, for any value of the momentum and the energy neighboring
$0$, there exists at least {\bfi  one RPO} with trivial isotropy.
An in--depth study of these motions  can be found in~\cite{cb}. 
\item [(ii)] {\bf The spring pendulum in three dimensions:} the
spring pendulum in space is a bob of mass $m$ under the action of a
constant gravitational field of acceleration $g$ attached to a
spring whose length at rest equals $l$.
The Hamiltonian function associated to this system is
\[
h(x,y,z,p_x,p_y, p _z)=\frac{1}{2}(p _x^2+ p _y^2+ p _z ^2)- m g
z+\sigma(x ^2+ y ^2 + (z - l) ^2,p _x^2+ p _y^2,z,p _z).
\]
where $ \sigma $ is a smooth real valued function such that
$\sigma(x ^2+ y ^2 + (z - l) ^2,p _x^2+ p _y^2,z,p _z)=\frac{k}{2}(x
^2+ y ^2 + (z - l) ^2) + $ terms of order higher or equal
than two in the first two variables and strictly higher than two in
$z$ and $p _z$. This system presents the same
$S ^1
$--symmetry as the previous example. 
The point $(0,0,l+
\frac{gm}{k},0,0,0)$ is a stable equilibrium where the
linearization of the Hamiltonian vector field has two pairs of
triple purely imaginary eigenvalues equal to $\pm
\frac{i}{2}\sqrt{\frac{k}{m}}$. Hence, in this case, the resonance
space is six dimensional. When there are no higher order terms in
the expansion of the function $\sigma$ we have a linear spring
(which does not mean that its associated Hamiltonian vector field is
linear) obeying Hooke's law with elastic constant
$k$. In this case, the spring pendulum is too degenerate and
does not satisfy hypothesis {\bf (H2)}.
However, in the presence of generic non linear terms in the
expansion of $\sigma$,  Theorem~\ref{theorem 1}
predicts the existence of {\bfi  two RPOs} with
trivial isotropy, for any value of the momentum and the energy
neighboring the values at the equilibrium.
\item [(iii)] {\bf An example with spherical symmetry:} consider $T
^\ast \Bbb R ^3\simeq \Bbb R^3\times \Bbb R^3$ with the canonical
symplectic form and the group ${\rm SO (3)}$ acting canonically by
diagonal transformations. If we denote by $({\bf q}, {\bf p} )$ the
elements of $\Bbb R ^3 \times \Bbb R^3 $, the momentum map
associated to this action is the standard angular momentum
$\mathbf{J}({\bf q}, {\bf p} )= {\bf q}\times {\bf p} $. Consider
the Hamiltonian function:
\[
h({\bf q}, {\bf p} )= a\| {\bf p}\| ^2 + b \| \mathbf{q}\| ^2+ f(\|
{\bf p}\| ^2, \| \mathbf{q}\| ^2, \mathbf{q}\cdot {\bf p}),
\]
where $a$ and $b$ are real strictly positive constants and $f$ is a
function whose Taylor expansion only contains terms of order two or
higher in its variables. The associated Hamiltonian vector field
has a stable equilibrium at the origin and its linearization at
that point exhibits one pair of triple purely imaginary eigenvalues
that therefore have a six dimensional resonance space associated.
We now look for RPOs nearby the equilibrium. The origin is
surrounded by points with isotropy subgroup $K$ either trivial or
equal to a circle. When $K= S ^1 $, the vectors $\mathbf{q}$ and
${\bf p} $ are forced to be parallel and therefore its momentum
value equal to zero, hence there is only room for planar periodic
solutions.
When  $K=\{e\} $, the vectors $\mathbf{q}$ and
${\bf p} $ are necessarily not parallel and therefore its
angular momentum is different from zero. Now, for any momentum
value  $\mu\in \Bbb R^3$ close to zero, the orbifold $
\mathbf{J}^{-1}(\mu)\cap Q _{\{e\}} ^{-1} (1)/(({\rm SO (3)})
_\mu\times  S ^1)$ is zero dimensional because $({\rm SO (3)})
_\mu\simeq S ^1 $ consequently, any function $f$ that depends
non--trivially on $\mathbf{q}\cdot {\bf p}$ makes the Hamiltonian
vector field associated to the corresponding $h$ exhibit, by
Theorem~\ref{theorem 1},  at least {\bfi  one RPO} with trivial
isotropy,  for any value of the momentum and the energy neighboring
the values at the equilibrium.
\quad $\blacklozenge $
\end{description}
\end{examples}

\begin{remark}
\normalfont
Note that in the absence of symmetries and for non trivial manifolds, that is $G=\{e\}$ and $\dim V>0$, 
the first part of the main estimate~(\ref{main estimate 1}) reduces
to 
\begin{equation}
\label{moser case}
\frac{1}{2}\left(\dim U\sno\right),
\end{equation}
which coincides with the conclusions of Moser's version~\cite{moser 76}  of Theorem~\ref{weinstein moser theorem}.\ \ \ $\blacklozenge$
\end{remark}

\begin{remark}
\normalfont 
The main feature of the first part of the 
estimate~(\ref{main estimate 1}), namely,
\[
\frac{1}{2}\left(\dim U\sno^K-\dim (N(K)/K)-\dim
\left(N(K)/K\right)_{\lambda}\right)
\]
is that it does not involve any dynamics, that is, neither the Hamiltonian $h$ nor any of its byproducts are present in it. The kinematical setup of the problem, in our case given by the symplectic representation of $G$ on $V$, fully determine the number of RPOs that can be expected around a stable symmetric equilibrium, induced by ANY Hamiltonian system that satisfies the
hypotheses of the Theorem. 

This dynamical independence, that was already present in the non symmetric result of Weinstein and Moser~(\ref{moser case}), has a price in terms of sharpness, since in general, the Morse theoretical part of~(\ref{main estimate 1})
involving the Euler characteristic is expected to  give better
results.   
\ \ \ $\blacklozenge$
\end{remark}

\begin{remark}
\normalfont
Despite the one half in front of the first part of estimate~(\ref{main estimate 1}) we always obtain an integer out of it. Indeed,  $U\sno^K$ is a symplectic vector subspace of the symplectic vector space $U\sno$, and therefore of even dimension. Also, the coadjoint orbit 
$(\elek)\cdot \lambda$ is also a symplectic manifold of even
dimension equal to $\dim
(N(K)/K)-\dim\left(N(K)/K\right)_{\lambda}$, hence $-\dim
(N(K)/K)-\dim\left(N(K)/K\right)_{\lambda}$ is necessarily also an
even number.\ \ \ $\blacklozenge$
\end{remark}

\subsection{Proof of Theorem~\ref{theorem 1}}

The proof of Theorem~\ref{theorem 1} is unfortunately quite long
and technical. In order to make it more accessible to the reader I
have divided it into several subsections that accomplish a
self--contained task.

The general strategy of the proof consists of  characterizing the 
relative periodic
orbits that we are looking for,
through the use of the normal form theorem,  as the critical points
of an equivariant function on certain manifold. The orbit space of
this manifold with respect to the natural group action that leaves
it invariant is a symplectic orbifold (the group action is locally
free) which allows us to implement the various critical point
theory techniques that we introduced in Section~\ref{Results on critical point theory of functions on 
compact manifolds}.

\subsubsection{The manifold $\J\selehdk\inv(\lambda)\cap Q\inv(1) $} 

We start the proof of the theorem by showing that for momentum
values $\lambda \in  \mathfrak{l}^\ast$ nearby zero, the manifolds
$\J\selehdk\inv(\lambda) $ and $Q\inv(1) $ intersect transversely
and, therefore, the intersection forms a manifold.  This manifold is
important because it will constitute the numerator of the
symplectic orbifold that we described in the previous paragraph.

\begin{proposition}
\label{proposition crucial}
Suppose that we are under the hypotheses of Theorem~\ref{theorem 1}. Let $K\subset G$ be an isotropy subgroup of the $G$--action on $U\sno$ and $\J\selehdk:(U\sno)_K\rightarrow\geles$ be the momentum map corresponding to the free $\elek$--action on $(U\sno)_K$. Suppose that the set $\J\selehdk\inv(0)\cap Q\shk\inv(1)$ is not empty. Then, there is a neighborhood $B(K)\subset\geles$ of $0$ in $\geles$ such that for any $\lambda\in B(K)$:
\begin{description}
\item[(i)] $\J\selehdk\inv(\lambda)\pitchfork Q\inv(1)=\J\selehdk\inv(\lambda)\pitchfork Q\shk\inv(1)$, that is, $\J\selehdk\inv(\lambda)$  intersect transversely with  $Q\inv(1)$ and with $Q\shk\inv(1)$.
\item[(ii)] $\J\selehdk\inv(\lambda)\cap Q\inv(1)=\J\selehdk\inv(\lambda)\cap Q\shk\inv(1)$ is a compact submanifold of $(U\sno)_K$  of dimension $\dim(U\sno)_K-\dim\elek-1=\dim U\sno\suhk-\dim\elek-1$.
\item[(iii)] The submanifold $\J\selehdk\inv(\lambda)\cap Q\inv(1)=\J\selehdk\inv(\lambda)\cap Q\shk\inv(1)\subset (U\sno)_K$ does not contain any $\elek$--relative equilibrium of the system $((U\sno)_K,\omega|_{(U\sno)_K},Q\shk)$.
\end{description}
\end{proposition}

\noindent\textbf{Proof}\ \ \textbf{(i)} Given that the 
$\elek$--action on $(U\sno)_K$ is free, the momentum map $\J\selehdk$
is a submersion and therefore its level sets are always
submanifolds of $(U\sno)_K$, and consequently of $U\sno$. At the
same time, the definiteness of  the quadratic form $Q$ implies that
its level sets are compact submanifolds of $U\sno$. We will show
the transversality of the level sets in the statement by showing
that  $\J\selehdk\inv(0)\pitchfork Q\inv(1)$. Since the
transversality is an open condition (see for instance the Stability
Theorem in page 35 of~\cite{guillemin pollack}), the result will
follow for $\J\selehdk\inv(\lambda)\pitchfork Q\inv(1)$, for
$\lambda$ close enough to zero. In order to prove that
$\J\selehdk\inv(0)\pitchfork Q\inv(1)$ we need to show that for each
$v\in\J\selehdk\inv(0)\cap Q\inv(1)$, which by hypothesis is not
empty, $T_v(U\sno)_K=T_v\J\selehdk\inv(0)+ T_v Q\inv(1)$. Since $T_v
Q\inv(1)$ has codimension one, it suffices to find a vector $w\in
T_v\J\selehdk\inv(0)$ that does not lie in $T_v Q\inv(1)$. This job
is done by  
$w=\ddto v+tv\in T_v(U\sno)\shk$ since a straightforward calculation
shows that, in one hand, $T_v\J\selehdk\cdot
w=0 $
because $v\in\J\selehdk\inv(0)$ and, at the same time,
$
\bd Q(v)\cdot w=\bd^2h(0)(v,v)=2Q\suhk(v)\neq 0,
$
by the definiteness hypotheses on $Q\suhk$, hence $w\notin T_v Q\inv(1)$.

\medskip

\noindent\textbf{(ii)} Since $\J\selehdk\inv(\lambda)\cap
Q\inv(1)=\J\selehdk\inv(\lambda)\cap Q\shk\inv(1)\subset
(U\sno)\shk$ is a smooth submanifold of $U\sno$ and $(U\sno)\shk$
is a smooth submanifold of $U\sno$, then
$\J\selehdk\inv(\lambda)\cap Q\inv(1)$ is a smooth submanifold of
$(U\sno)\shk$.  We show the compactness of
$\J\selehdk\inv(\lambda)\cap Q\inv(1)$ by proving that it is
sequentially compact, that is, any sequence $\{x_n\}\subset
\J\selehdk\inv(\lambda)\cap Q\inv(1)$ has a convergent
subsequence in $\J\selehdk\inv(\lambda)\cap Q\inv(1)$. Indeed,
since $Q\inv(1)$ is a compact subset of $U\sno$, there exists a
subsequence $\{x_{n_k}\}$ of $\{x_n\}$, convergent in $Q\inv(1)$,
that is $x_{n_k}\rightarrow l$, with $l\in Q\inv(1)$. Since
$\{x_{n_k}\}\subset \J\selehdk\inv(\lambda)$ for all $n_k$ and
$\J\selehdk\inv(\lambda)$ is closed then $l\in
\J\selehdk\inv(\lambda)$ necessarily, as required. The
claim on the dimension is a consequence of the Transversality
Theorem.

\medskip

\noindent\textbf{(iii)} We proceed by contradiction: suppose that the point $v\in (U\sno)_K$ is a relative equilibrium of the system with velocity $\xi\in\gele$. This implies that $\bd\left(Q\shk-\J\selehdk\suxi\right)(v)=0$. The Lagrange Multipliers Theorem (taking $\xi$ in Proposition~\ref{lagrange multipliers} as the Lagrange multiplier) implies that $\left.\bd Q\shk(v)\right|_{\ker T_v\J\selehdk}=0$. However, the transversality $\J\selehdk\inv(\lambda)\pitchfork Q\shk\inv(1)$ implies the existence of vectors $w\in \ker T_v\J\selehdk$ 
that do not belong to $T_v Q\shk\inv(1)$ and therefore $\bd
Q\shk(v)\cdot w\neq 0$, which represents a contradiction.\ \ \
$\blacksquare$

\medskip

\subsubsection{Normal form formulation}

As we already said, the main idea behind the
theorem consists of using the normal form
theorem and Proposition~\ref{elementary result} to reduce the search
for RPOs with isotropy subgroup $K$ to the search for periodic
orbits of the Hamiltonian systems of the form
$(\omega|_{V\shk},h_\xi)$ with $h_\xi:=h|_{V\shk}-\J\suxi\selehdk$. We
will first apply the ideas introduced in Section~\ref{The resonance
space and normal form reduction} to construct a normal form for
these systems that will set up the problem of this search for periodic orbits
in terms of the search for relative equilibria of a $S^1$--action
that we will describe in what follows.

Firstly, consider the symmetric Hamiltonian system $(V\suhk,\omega|_{V\suhk},h|_{V\suhk}, N(K)/K,\J\selehuk:V\suhk\rightarrow\geles)$. The RPOs of this system amount to RPOs of the original system whose isotropy subgroups include $K$. Let 
$A\suhk:=DX_{h|_{V\suhk}}(0)$ be the linearization at zero of the
Hamiltonian vector field $X_{h|_{V\suhk}}$. By the hypotheses on
$X_h$, the eigenvalues $\pm i\no$ are in the spectrum of $A\suhk$,
and the corresponding resonance space is $U\sno\suhk$. The mapping
$\theta\in S^1\mapsto{\rm
e}^{\frac{\theta}{\no}A\suhk}|_{U\sno\suhk}$ generates a symplectic
linear $S^1$--action on  $(U\sno\suhk,\omega|_{U\sno\suhk})$, whose
momentum map is given by $\frac{1}{\no}Q\suhk$.

Consider now $h_\xi\suhk:=h|_{U\sno\suhk}-\J\selehuk\suxi$, with
$\xi\in\gele$,  as a $\gele$--parameter family of Hamiltonian
functions on $U\sno\suhk$ in the sense of point~\textbf{(v)} in
Section~\ref{The resonance space and normal form reduction}. Since
for any $\xi\in\gele$ this family satisfies that $h_\xi\suhk(0)=0$,
$\bd h_\xi\suhk(0)=0$, and
$DX_{h_0\suhk}(0)=DX_{h|_{U\sno^K}}(0)=A|_{U\sno\suhk}=A\suhk$ is non
degenerate, we can construct a normal form equivalent system  
$\widehat{h_\xi\suhk}$ whose $S^1$--relative equilibria give us the
periodic orbits of $h_\xi\suhk$. Due to the fact that the $\elek$ and
$S^1$--actions on $U\sno\suhk$ commute and that $\J\selehuk$ is
quadratic, the normal form $\widehat{h_\xi\suhk}$ can be chosen so
that
\[
\widehat{h_\xi\suhk}=\widehat{h|_{U\sno\suhk}}-\J\selehuk\suxi,
\]
with $\widehat{h|_{U\sno\suhk}}$ an $\elek\times S^1$--invariant function on 
$U\sno\suhk$ such that 
\begin{equation}
\label{higher}
\widehat{h|_{U\sno\suhk}}(u)=Q^K(u)+\frac{1}{2}\bd^2 h|_{U\sno\suhk}(0)(u,u)+\cdots+\frac{1}{k!}\bd^k h|_{U\sno\suhk}(0)(u^{(k)})+\text{(higher order terms)}. 
\end{equation}

\subsubsection{The critical points of the normal form Hamiltonian
and a blow up argument}

In the following lemma we evaluate the critical points of the
restriction of the function $\widehat{h|_{U\sno\suhk}}$ to the level
sets of the form $\J\selehdk\inv(\lambda)\cap Q_K\inv(\epsilon)$,
where $\lambda\in B(K)$, the neighborhood of zero in $\geles$
introduced in Proposition~\ref{proposition crucial}, and
$\epsilon>0$ is very small. Furthermore, a blow up argument will
show that  these critical points can be arranged in smooth branches.
\begin{lemma}
\label{smooth branches}
Suppose that we are under the hypotheses of Theorem~\ref{theorem 1}. Then, the restriction of the function $h_k\in C\suinf(U\sno\suhk)$ defined by $h_k(u):=\frac{1}{k!}\bd^k h|_{U\sno\suhk}(0)(u^{(k)})$,  to the level sets of the form $\J\selehdk\inv(\lambda)\cap Q_K\inv(1)$, where $\lambda\in B(K)$, the neighborhood of zero in $\geles$ introduced in Proposition~\ref{proposition crucial}, has at least
\begin{equation}
\label{critical lusternik}
{\rm max}\left[{\rm Cat}\left(\J\selehdk\inv(\lambda)\cap Q_K\inv(1)/(L_\lambda\times S^1)\right), \chi\left(\J\selehdk\inv(\lambda)\cap Q\shk\inv(1)\right)^{L_\lambda\times S^1}\right]
\end{equation}
distinct critical points. 

Furthermore, let $\lo\in B(K)\subset\geles$ be arbitrary but fixed and let $\{u_1,\ldots,u_k\}$ be the set of critical points of the restriction of the function $h_k$ to the level set $\J\selehdk\inv(\lo)\cap Q_K\inv(1)$, provided by~(\ref{critical lusternik}). Then, for each $u_i$, $i\in\{1,\ldots,k\}$, there exist a neighborhood $E\subset\mathbb{R}$ of the origin in the real line, a neighborhood $V\slo\subset\geles$ of $\lo$ in $\geles$, and a smooth function $\rho:E\times V\slo\rightarrow Q\shk\inv(1)$ such that $\rho(0,\lo)=u_i$ and also, the function $v:E\times V\slo\rightarrow U\sno\suhk$ defined by $v(r,\lambda):=r\rho(r,\lambda)$ satisfies that:
\begin{enumerate}
\item[{\bf (i)}] $v(r,\lambda)\in (U\sno)\shk$ iff $r\neq 0$.
\item[{\bf (ii)}] $Q\suhk(v(r,\lambda))=r^2$ and $\J\selehuk(v(r,\lambda))=r^2\lambda$, $r\in E$, $\lambda\in V\slo$.
\item[{\bf (iii)}] $\bd \widehat{h|_{U\sno\suhk}}|_{\J\selehuk\inv(r^2\lambda)\cap (Q^K)\inv(r^2)}(v(r,\lambda))=0$, that is, the branches $v(r,\lambda)$ are made of critical points of the restriction of the function $\widehat{h|_{U\sno\suhk}}$ to the level sets $\J\selehuk\inv(r^2\lambda)\cap (Q^K)\inv(r^2)$.
\end{enumerate}
\end{lemma}

\noindent\textbf{Proof} Firstly, note that the estimate~(\ref{critical lusternik}) is a straightforward consequence of Proposition~\ref{v manifolds},~(\ref{euler estimate}), and the invariance properties of $h_k$.

We now prove the existence of the branches in the statement. Let $u\szero$ be one of the critical points of the restriction of the function $h_k$ to the level set $\J\selehdk\inv(\lo)\cap Q_K\inv(1)$, provided by~(\ref{critical lusternik}). The transversality of $\J\selehdk$ and $Q\shk$ implies the existence of a very convenient coordinate patch around $u\szero$ in $Q\shk\inv(1)$. In order to construct it, let us show first that the restriction $\J\selehdk|_{Q\shk\inv(1)\cap\J\selehdk\inv(B(K))}$ is a submersion onto its image. Indeed, for any $u\in Q\shk\inv(1)\cap\J\selehdk
\inv(B(K))$ we have that
\begin{align*}
\dim{\rm Im}&\left(T_u\J\selehdk|_{Q\shk\inv(1)\cap\J\selehdk\inv(B(K))}\right)=\dim\left(\ker T_u Q\shk\right)-\dim\left(\ker T_u\J\selehdk|_{Q\shk\inv(1)}\right)\\
	&=\dim U\sno\suhk-1-\dim\left(\ker T_u\J\selehdk\cap\ker T_u Q\shk\right)\\
	&=\dim U\sno\suhk-1+\dim\left(\ker T_u\J\selehdk+\ker T_u Q\shk\right)-\dim\left(\ker T_u\J\selehdk\right)-\dim\left(\ker T_u Q\shk\right)\\
	&=\dim\geles,
\end{align*}
as required. In the last equality we used that since $u\in \J\selehdk\inv(B(K))$, then $\ker T_u\J\selehdk\pitchfork\ker T_u Q\shk$, and therefore $\dim\left(\ker T_u\J\selehdk+\ker T_u Q\shk\right)=\dim U\sno\suhk$. In these circumstances, the Local Submersion Theorem (see for instance~\cite[Theorem 3.5.2]{mta}) implies the existence of a neighborhood $V\slo$ of $\lo$ in $\geles$, a neighborhood $W$ of the origin in $\mathbb{R}^s$, with $s=\dim U\sno\suhk-\dim\geles-1$, and a mapping $\varphi:V\slo\times W\rightarrow Q\shk\inv(1)$ that is a diffeomorphism onto its image, such that:
\begin{equation}
\label{properties phi}
\varphi(\lo,0)=u\szero\quad\text{and}\quad\J\selehdk(\varphi(\lambda,w))=\lambda,\qquad \lambda\in V\slo, w\in W.
\end{equation}
We will further improve this coordinate patch around $u\szero$ by "factoring out" the $L\slo\times S^1$--action in it. Indeed, since the Lie group $L\slo\times S^1$ acts on $Q\shk\inv(1)$, we can induce a local action of this group on $V\slo\times W$ by declaring $\varphi$ to be equivariant, that is, for any $g=(l,\theta)\in  L\slo\times S^1$ close enough to the identity we define
$
g\cdot (\lambda,w):=\varphi\inv\left(g\cdot\varphi(\lambda,w)\right).
$
Note that by the very definition of this action and by~(\ref{properties phi}), we have that
\begin{equation}
\label{j and g}
\J\selehdk(\varphi(g\cdot(\lambda,w)))=\J\selehdk(\varphi(\varphi\inv(g\cdot\varphi((\lambda,w)))))=\J\selehdk(g\cdot\varphi(\lambda,w))=l\cdot \J\selehdk(\varphi(\lambda,w))=l\cdot\lambda,
\end{equation}
consequently 
\begin{equation}
\label{local action}
g\cdot(\lambda,w)=(l,\theta)\cdot(\lambda,w)=(l\cdot\lambda, \Phi(g,\lambda,w)),
\end{equation} 
for some smooth function $\Phi$. The dot in $l\cdot\lambda$ denotes 
the coadjoint action of $L\slo$ on $\geles$. Note that in~(\ref{j
and g}) we used the equivariance of $\J\selehdk$ with respect to the
$L\slo$--action that was inherited from its $L$--equivariance, as
well as the invariance of $\J\selehdk$ with respect to the
$S^1$--action which is justified, via Noether's Theorem, by the fact
that this action is induced by a Hamiltonian flow associated to a
$L$--invariant Hamiltonian, namely
$Q\suhk$.  

Using these remarks we are now going to construct a slice for the local $L\slo\times S^1$--action on $V\slo\times W$ that goes through the point $(\lo,0)$. Firstly, it is easy to see by using~(\ref{local action}) that
\[
T_{(\lo,0)}\left((L\slo\times S^1)\cdot(\lo,0)\right)=\{0\}\times W',
\]  
where $W'\subset W$ is a vector subspace of $W$. The remarks that we made in Section~\ref{Proper actions, fixed 
point sets, slices, and normalizers.} about the construction of the slices implies the existence of a smooth mapping $\psi$ diffeomorphic onto its image of the form
\[
\begin{array}{cccc}
\psi:&V\slo\times U&\longrightarrow&V\slo\times W\\
	&(\lambda,u)&\longmapsto&(\lambda,\eta(\lambda, u))
\end{array}
\]
whose image (shrinking $V\slo$ if necessary) is a local slice through $(\lo,0)$ for the local $L\slo\times S^1$--action on $V\slo\times W$. The set $U$ is an open neighborhood of the origin in a vector space isomorphic to $W/W'$ and $\eta:V\slo\times U\rightarrow V\slo\times W$ is a smooth map such that $\eta(\lo,0)=0$.

In these circumstances, the version of the Slice Theorem that we presented in Section~\ref{Proper actions, fixed 
point sets, slices, and normalizers.} implies the existence of a local cross--section $\sigma:\mathcal{Z}\subset L\slo\times S^1/(L\slo\times S^1)_{u\szero}\rightarrow L\slo\times S^1$ of the homogeneous space $L\slo\times S^1/(L\slo\times S^1)_{u\szero}$, and a smooth map $F$ of the form 
\begin{equation}
\label{slice patch}
\begin{array}{cccc}
F:&\mathcal{Z}\times V\slo\times U&\longrightarrow&Q_K\inv(1)\\
	&(z,\lambda,u)&\longmapsto&\sigma(z)\cdot\varphi(\lambda,\eta(\lambda, u))
\end{array}
\end{equation}
that is a diffeomorphism onto an open set of  $Q_K\inv(1)$ that contains $u\szero$.

We will now use this coordinate patch to obtain the branches $v(r,\lambda)$ whose existence we claim in the second part of the statement. We start by setting up the problem in polar coordinates since it is a polar blowing--up argument what will give us the result. Let $2n$ be the dimension of the symplectic vector space $U\sno\suhk$. We will denote by $S^{2n-1}$ the sphere in $U\sno\suhk$ obtained by using the norm associated to the definite quadratic form $Q\suhk$. We now define the blown--up Hamiltonian 
$h ^b:\mathbb{R}\times S^{2n-1}\rightarrow\mathbb{R}$, as:
\[
h ^b(r,u):=\widehat{h|_{U\sno\suhk}}(ru),\qquad r\in\mathbb{R}, \quad u\in S^{2n-1}.
\]
When the variable $u$ is in a neighborhood of $u\szero$, we can use~(\ref{slice patch}) to give a local expression for $h ^b$ in $\mathcal{Z}\times V\slo\times U$--variables; 
let $h ^{lb}:\mathbb{R}\times \mathcal{Z}\times V\slo\times
U\rightarrow\mathbb{R}$ be the local expression of the blown--up
Hamiltonian, defined by
\[
h ^{lb}(r,z,\lambda,u):=h ^b(r(\sigma(z)\cdot\varphi(\lambda,\eta(\lambda, u))))=\widehat{h|_{U\sno\suhk}}(r(\sigma(z)\cdot\varphi(\lambda,\eta(\lambda, u)))).
\]
Notice that the $G$--invariance of the Hamiltonian $h$ implies that
\[
h ^{lb}(r,z,\lambda,u)=\widehat{h|_{U\sno\suhk}}(r(\sigma(z)\cdot\varphi(\lambda,\eta(\lambda,
u))))=\widehat{h|_{U\sno\suhk}}(r\varphi(\lambda,\eta(\lambda,
u)))\equiv h ^{lb}(r,\lambda,u),
\]
that is, $h ^{lb}$ does not depend on the $\mathcal{Z}$--variables. The main advantage of the use of these coordinates is the fact we can search for the critical points of the restriction of $\widehat{h|_{U\sno\suhk}}$ to the level sets $\J\selehuk\inv(\lambda)\cap (Q^K)\inv(\epsilon)$ by looking for the triples $(r,\lambda,u)\in \mathbb{R}\times V\slo\times U$ for which 
$D_Uh ^{lb}(r,\lambda,u)=0$ ($D_Uh ^{lb}$ denotes the partial Fr\'echet
derivative of $h ^{lb}$ relative to the $U$--variable). More specifically: 
$D_Uh ^{lb}(r,\lambda,u)=0$ iff the restriction of $\widehat{h|_{U\sno\suhk}}$ to the level set
$\J\selehuk\inv(r^2\lambda)\cap (Q^K)\inv(r)$ has a critical point at $r\varphi(\lambda,\eta(\lambda, u))$. Using now hypotheses~\textbf{(H1)} and~\textbf{(H2)} on the Hamiltonian $h$ we can write
\begin{equation}
\label{first equivalence}
h ^{lb}(r,\lambda,u)=r^2+\alpha_4r^4+\alpha_6r^6+
\cdots+\alpha_{k-2}r^{k-2}+r^k h_k^{lb}(\lambda,u)+o(r^k,\lambda,u),
\end{equation} 
where $\alpha_4,\alpha_6,\ldots,\alpha_{k-2}$ are real 
coefficients and
$h_k^{lb}(\lambda,u):=h_k(\varphi(\lambda,\eta(\lambda, u)))$.
Expression~(\ref{first equivalence}) can be rewritten of the form
\[
h ^{lb}(r,\lambda,u)-r^2-\alpha_4r^4-\alpha_6r^6-\cdots-\alpha_{k-2}r^{k-2}=r^kg(r,\lambda,u),
\]
with $g$ is a smooth function on his variables such that 
\begin{equation}
\label{reduction to g}
D_Ug(r,\lambda,u)=0\qquad\text{if and only if}\qquad D_Uh ^{lb}(r,\lambda,u)=0.
\end{equation}
The Taylor expansion of $g$ on the $r$ variables around $r=0$ has the form
\[
g(r,\lambda,u)=h_k^{lb}(\lambda,u)+o(r,\lambda,u)=h_k(\varphi(\lambda,\eta(\lambda,
u)))+o(r,\lambda,u)
\]
Notice that
\[
g(0,\lo,0)=h_k(\varphi(\lo,\eta(\lo, 0)))=h_k(u\szero),
\]
and that 
$$D_Ug(0,\lo,0)=D_U h _k^{lb}(\lo,0)=\bd h_k(u\szero)\cdot(D_W\varphi(\lo,0)\cdot D_U\eta(\lo,0))=0,$$ since $D_W\varphi(\lo,0)\cdot D_U\eta(\lo,0)$ maps into $T_{u\szero}(\J\selehdk\inv(\lo)\cap Q_K\inv(1))$, which, by the choice of $u_0$ lies in the kernel of $\bd h_k(u\szero)$ (recall that $u\szero$ was chosen to be a critical point of the restriction of $h_k$ to $\J\selehdk\inv(\lo)\cap Q_K\inv(1)$). Moreover, it is easy to see that for any pair $u,v\in U$:
\[
D^2_Ug(0,\lo,0)(u,v)=\bd^2 h_k|_{\J\selehdk\inv(\lo)\cap Q_K\inv(1)}(u_0)(D_W\varphi(\lo,0)\cdot D_U\eta(\lo,0)\cdot u,D_W\varphi(\lo,0)\cdot D_U\eta(\lo,0)\cdot v).
\]
The $(N (K)/K)_\lambda\times S ^1$--Morse condition on
$h_k|_{\J\selehdk\inv(\lo)\cap Q_K\inv(1)}$ implies that
$D^2_Ug(0,\lo,0)$ is a non degenerate quadratic form since the
image of the linear mapping
$D_W\varphi(\lo,0)\cdot D_U\eta(\lo,0):U\rightarrow
T_{u\szero}(\J\selehdk\inv(\lo)\cap Q_K\inv(1))$ is a vector
subspace of
$T_{u\szero}(\J\selehdk\inv(\lo)\cap Q_K\inv(1))$ that, by
construction, is complementary to $T_{u\szero}((L\slo\times
S^1)\cdot u\szero)$. In this situation the Implicit Function
Theorem guarantees the existence of a function $u:E\times
V\slo\rightarrow U$ (shrink $V\slo$ if necessary) with
$E\subset\mathbb{R}$ a neighborhood of the origin in $\mathbb{R}$,
such that 
\[
D_Ug(r,\lambda,u(r,\lambda))=0,
\] 
which, by~(\ref{reduction to g}) is equivalent to having $D_Uh ^{lb}(r,\lambda,u(r,\lambda))=0$ and consequently
\[
\bd\widehat{h|_{U\sno\suhk}}|_{\J\selehuk\inv(r^2\lambda)\cap (Q^K)\inv(r)}(r\varphi(\lambda,\eta(\lambda, u(r,\lambda))))=0.
\]
The claim of the Lemma follows by taking $\rho(r,\lambda):=\varphi(\lambda,\eta(\lambda, u(r,\lambda)))$ and $v(r,\lambda):=r\rho(r,\lambda)=r\varphi(\lambda,
\eta(\lambda, u(r,\lambda)))$. \quad $\blacksquare$

\subsubsection{Lagrange multipliers and RPOs}

We now show that the branches of critical points found in
Lemma~\ref{smooth branches} amount to RPOs of the original system. 

Let now $\lo\in\gele$ and $v(r,\lambda)$ (with $v:E\times V\slo\rightarrow U\sno^K$) be one of the
branches of critical points of the restriction of the function $\widehat{h|_{U\sno\suhk}}$ to the level sets $\J\selehuk\inv(r^2\lambda)\cap (Q^K)\inv(r^2)$ provided by~(\ref{critical lusternik}).
Proposition~\ref{lagrange multipliers} guarantees, for each $v(r,\lambda)$ the existence of a multiplier $(\Lambda(r,\lambda),c(r,\lambda))\in\gele\times\mathbb{R}$ such that:
\begin{equation}
\label{identity pairing}
\bd \widehat{h|_{U\sno\suhk}}(v(r,\lambda))=c(r,\lambda)\bd Q^K(v(r,\lambda))+\bd\J\selehuk^{\Lambda(r,\lambda)}(v(r,\lambda)), 
\end{equation}
which implies that the point $v(r,\lambda)$ is a periodic point of $X_{\widehat{h\suhk_{\Lambda(r,\lambda)}}}$. If we are able to show that $\Lambda(r,\lambda)$ can be made very small so that we can use the Normal Form Theorem, all these periodic points will amount to periodic orbits of $X_{h|_{U\sno\suhk}-\J\selehuk^{\Lambda(r,\lambda)}}$ and, by Proposition~\ref{elementary result}, to RPOs of $X_h$. We will prove this point by pairing both sides of~(\ref{identity pairing}) at the point $v(r,\lo)$ with the vector $w_\lambda:=D_\lambda v(r,\lo)\cdot\lambda$, $\lambda\in\geles$ arbitrary, and taking into account that by the very construction of the function $v(r,\lambda)$, $\J\selehuk(v(r,\lambda))=r^2\lambda$, for any $r\in E$ and any $\lambda\in V\slo$. Indeed, 
\begin{equation}
\label{almost there}
\bd \widehat{h|_{U\sno\suhk}}(v(r,\lo))\cdot w_\lambda=c(r,\lo)\bd Q^K(v(r,\lo))\cdot w_\lambda+\bd\J\selehuk^{\Lambda(r,\lo)}(v(r,\lo))\cdot w_\lambda.
\end{equation}
Note first that $\bd Q^K(v(r,\lo))\cdot w_\lambda=0$ for all $\lambda\in\geles$. Also,
\begin{multline}
\label{almost there 1}
\bd\J\selehuk^{\Lambda(r,\lo)}(v(r,\lo))\cdot w_\lambda=\ddto\J\selehuk^{\Lambda(r,\lo)}(v(r,\lo+t\lambda))\\
=\ddto\langle r^2(\lo+t\lambda),\Lambda(r,\lo)\rangle=r^2\langle\lambda,\Lambda(r,\lo)\rangle.
\end{multline}
This equality implies, together with hypotheses~\textbf{(H1)} and~\textbf{(H2)} that the multiplier $\Lambda(r,\lambda)$ is a smooth function in its variables $\Lambda:E\times V\slo\rightarrow\gele$ since it can be written as a composition of smooth functions, namely
\begin{equation}
\label{not really}
\Lambda(r,\lambda)=\frac{1}{r^2} \bd \widehat{h|_{U\sno\suhk}}(v(r,\lambda))\cdot D_\lambda v(r,\lambda).
\end{equation}
Even though in the previous equality it seems that there is a singularity at $r=0$ we see in what follows that it is not the case. Indeed,  by hypothesis~\textbf{(H2)} we can write
\begin{equation}
\label{useful}
\widehat{h|_{U\sno\suhk}}(v(r,\lambda))=r^2+\alpha_4r^4+\alpha_6r^6+\cdots+\alpha_{k-2}r^{k-2}+r^k h_k(\rho(r,\lambda))+o(r^k,\lambda),
\end{equation}
where $\alpha_4,\alpha_6,\ldots,\alpha_{k-2}$ are real coefficients. It is easy to see from~(\ref{useful}) that for any $\eta\in\geles$
\[
\bd \widehat{h|_{U\sno\suhk}}(v(r,\lambda))\cdot w_\eta=r^k\bd h_k(\rho(r,\lambda))\cdot D_\lambda\rho(r,\lambda)\cdot\eta+o(r^k,\lambda).
\]
Given that by hypothesis~\textbf{(H2)} $k\geq 4$, expression~(\ref{not really}) can be rewritten as
\[
\Lambda(r,\lambda)=r^{k-2}\bd h_k(\rho(r,\lambda))\cdot D_\lambda\rho(r,\lambda)\cdot\eta+o(r^{k-2},\lambda),
\]
where the smoothness of the function $\Lambda(r,\lambda)$ is apparent as well as the fact that $\Lambda(0,\lambda)=0$. These two points together imply that the multiplier $\Lambda(r,\lambda)$ can be made as small as we want by taking $r$ sufficiently close to the origin, as desired. This allows us to use the Normal Form Theorem to conclude that for $r$ small enough, the RPOs of 
$\widehat{h|_{U\sno\suhk}}$ amount to RPOs of the original system.

We will now prove that the RPOs that we just obtained have relative periods close to $T\sno$ by showing that as $r$ tends to zero, the multiplier $c(r,\lambda)$ approaches to 1.
We prove this point by pairing both sides of~(\ref{identity pairing})  with the vector $u:=D_r v(r,\lambda)\cdot 1$. First of all
\begin{equation}
\label{yes 1}
\bd Q\suhk(v(r,\lambda))\cdot u=\ddto Q\suhk(v(r+t,\lambda))=\ddto (r+t)^2=2r.
\end{equation}
Also,
\begin{multline}
\label{yes 2}
\bd\J^{\Lambda(r,\lambda)}\selehuk(v(r,\lambda))\cdot u
=\ddto\langle\J\selehuk(v(r+t,\lambda)),\Lambda(r,\lambda)\rangle\\
=\ddto\langle(r+t)^2\lambda,\Lambda(r,\lambda)\rangle=2r\langle\lambda,\Lambda(r,\lambda)\rangle.
\end{multline}
From expression~(\ref{useful}) it is easy to see that
$\bd \widehat{h|_{U\sno\suh}}(v(r,\lambda))\cdot u=2r+o(r,\lambda)$. This equality, together with~(\ref{yes 1}) and~(\ref{yes 2}) imply, when substituted in~(\ref{identity pairing}) paired with $u$ that 
\[
1-c(r,\lambda)=\langle\lambda,\Lambda(r,\lambda)\rangle+o(r,\lambda).
\]
Since $\Lambda(r,\lambda)\rightarrow 0$ when $r$ tends to zero then $c(r,\lambda)\rightarrow 1$ as $r\rightarrow 0$, as claimed.

\subsubsection{Dimensional estimate of the Lusternik--Schnirelmann
category}

In order to conclude the proof of the theorem we will give a dimensional estimate of the first part of the estimate on the branches~(\ref{critical lusternik}) whose relatively periodic character we just proved. In order to obtain the claim~(\ref{main estimate 1}) in the statement of the theorem we just need to show that
\begin{equation}
\label{critical lusternik inequality}
{\rm Cat}\left(\J\selehdk\inv(\lambda)\cap Q_K\inv(1)/(L_\lambda\times S^1)\right)\geq\frac{1}{2}
\left[\dim(U\sno)_K-\dim(\elek)-\dim(\elek)_\lambda\right].
\end{equation}
We will prove this inequality with the help of 
Proposition~\ref{result b v manifolds}, taking in its statement
$\J\selehdk\inv(\lambda)\cap Q_K\inv(1)$ as the submanifold $N$ and
$L_\lambda\times S^1$ the subgroup $K$. Given that the mapping
$\J\selehdk\times \frac{1}{\no}Q_K$ is the momentum map
corresponding to the $\elek\times S^1$--action on $(U\sno)_K$, the
Reduction Lemma (see for instance~\cite[Lemma 4.3.2]{fom})
guarantees the technical hypotheses~(\ref{reduction hypotheses}),
namely:
\[
(T _z (\J\selehdk\inv(\lambda)\cap Q_K\inv(1))) ^\omega=(\ker T
_z(\J\selehdk\times \frac{1}{\no}Q_K))^\omega=T _z((\elek\times
S^1)\cdot z), \quad\text{as well as,}
\]
\[
T _z (\J\selehdk\inv(\lambda)\cap Q_K\inv(1))\cap T _z((\elek\times
S^1)\cdot z)=T _z((L _\lambda\times
S^1)\cdot z),
\]
for any $z \in \J\selehdk\inv(\lambda)\cap Q_K\inv(1) $.
In order to apply Proposition~\ref{result b v manifolds} we
need the isotropy subgroups of any point in
$\J\selehdk\inv(\lambda)\cap Q_K\inv(1)$ to be a finite subgroup of
$L_\lambda\times S^1$. This is so by the freeness and the local
freeness of the $L_\lambda$ and $S^1$--actions on $(U\sno)_K$,
respectively, and by the fact proved in the third part of
Proposition~\ref{proposition crucial} that 
$\J\selehdk\inv(\lambda)\cap Q_K\inv(1)$ does not contain any
$\elek$--relative equilibria of the system
$((U\sno)_K,\omega|_{(U\sno)_K},Q|_{(U\sno)_K})$ that would be the
only elements that could make the $L_\lambda\times S^1$--action not
be  locally free. Consequently, by Proposition~\ref{result b
v manifolds}, there is a  cohomology class $\theta \in
K ^2(\J\selehdk\inv(\lambda)\cap Q_K\inv(1)/(L_\lambda\times S^1);
\Bbb R)$ such that $\theta^k\neq 0$, with $k= \frac{1}{2}(\dim
(\J\selehdk\inv(\lambda)\cap Q_K\inv(1))-\dim (L_\lambda\times
S^1))$. Therefore, the cuplength of
$\J\selehdk\inv(\lambda)\cap Q_K\inv(1)/(L_\lambda\times S^1)$ is at
least
\[
\frac{1}{2}\left[\dim \left(\J\selehdk\inv(\lambda)\cap Q_K\inv(1)\right)-\dim L_\lambda-1\right]=\frac{1}{2}\left[\dim (U\sno)_K-\dim L-\dim L_\lambda-2\right],
\]
and hence, Proposition~\ref{cup length} establishes the
inequality~(\ref{critical lusternik inequality}), which concludes
the proof of the theorem.\ \ \ $\blacksquare$

\medskip

\begin{remark}
\label{rpos and not equilibria}
\normalfont
The choice of the neighborhood $B(K)$ in Proposition~\ref{proposition crucial} guarantees that the RPOs that we found in our theorem are nontrivial, that is, they are not just relative equilibria. Indeed, when $\lambda\in B(K)$, the level set of  $Q\shk\inv(1)$ intersects transversely the level set 
$\J\selehdk\inv(\lambda)$, which in part~\textbf{(iii)} of 
Proposition~\ref{proposition crucial} allowed us to prove the
inexistence of relative equilibria of the Hamiltonian vector field
$X_{Q\shk}$ with respect to the $\elek$--symmetry in
$\J\selehdk\inv(\lambda)\cap Q\shk\inv(1)$. Since the Hamiltonian
$h|_{(U\sno)\shk}$ can be considered as a perturbation of $Q\shk$ for
vectors of norm $r$ small enough, the transversality (which is a
stable property) of $h|_{(U\sno)\shk}\inv(r)$ with respect to
$\J\selehdk\inv(r^2\lambda)$ is still valid, and therefore the
inexistence of $\ele$--relative equilibria in
$\J\selehdk\inv(r^2\lambda)\cap h|_{(U\sno)\sh}\inv(r)$ as well.\ \
\ $\blacklozenge$
\end{remark}

\begin{remark}
\normalfont
Even though the hypothesis~\textbf{(H1)} appears in the proof of 
the theorem as a technical necessity, it turns out when it fails,
the existence of genuine RPOs that are not
either relative equilibria or plain periodic orbits is not possible.
Indeed, suppose that the restriction $h|_{U\sno^K}$ is purely
radial. In that case, $h|_{U\sno^K}$ is directly in normal form and
there exists a real smooth function
$f:\mathbb{R}\rightarrow\mathbb{R}$ such that
$h|_{U\sno^K}(v)=f(Q^K(v)),\,v\in U\sno^K$. In these circumstances,
expression~(\ref{identity pairing}) reduces to 
\[
f'(Q^K(v(r,\lambda)))\bd Q^K(v(r,\lambda))=c(r,\lambda)\bd Q^K(v(r,\lambda))+\bd\J\selehuk^{\Lambda(r,\lambda)}(v(r,\lambda)), 
\]
which amounts to 
\[
\bd h|_{U\sno^K}(v(r,\lambda))-\bd\J\selehuk^{\frac{\Lambda(r,\lambda)f'(Q^K(v(r,\lambda)))}{f'(Q^K(v(r,\lambda)))-c(r,\lambda)}}(v(r,\lambda))=0, 
\]
that is, in the absence of  hypothesis~\textbf{(H1)} $v(r,\lambda)$ is a branch of relative equilibria of the Hamiltonian vector field associated to $h$ (the reader interested in the technology for searching relative equilibria in the hypotheses of Theorem~\ref{theorem 1}, or even weaker, can check with~\cite{pascal, singular moser}, and references therein). Note that a trivial corollary of this comment is  that if $h$ is just quadratic and therefore its associated Hamiltonian vector field is linear, then there are no genuine RPOs associated to its dynamics. On other words, \emph{the relative periodic orbits around stable equilibria are purely non linear phenomena}.\ \ \ $\blacklozenge$ 
\end{remark}

\subsection{Relative periodic orbits with prescribed spatiotemporal
symmetry}

In the statement of Theorem~\ref{theorem 1} we optimized the search for the RPOs of our system by looking for them within the fixed point spaces corresponding to the isotropy subspaces of the $G$--action on $U\sno$. In Section~\ref{The resonance space and normal form reduction} we showed that the resonance space $U\sno$ can actually be endowed with a $G\times S^1$--symmetry which obviously contains more isotropy subgroups than merely the $G$--symmetry and that we could therefore utilize to obtain additional relatively periodic solutions. However, the reader should not forget that the $S^1$--symmetry is a feature owned solely by the system in normal form; the real system that we are dealing with is not $S^1$--symmetric. This fact does not pose a problem since the morphism that relates the $S^1$--relative equilibria of the normal form equivalent system to the periodic orbits of the original system transforms the $S^1$--symmetry of the normal form Hamiltonian into a $S^1$--symmetry \emph{of the periodic solutions of the original system}. The $G\times S^1$--action on the $T\sno$--periodic solutions is defined as $(g,\theta)\cdot u(t):=g\cdot u(t+\frac{t\theta}{2\pi})$, where $u:\mathbb{R}\rightarrow V$ is a smooth function such that $u(t+T\sno)=u(t)$.

The use of the isotropy subgroups of the $G\times S^1$--symmetry of the normal form equivalent system has been very fruitful in the symmetric bifurcation theory (see~\cite{g2} for a taste of it).

In this section we will generalize Theorem~\ref{theorem 1} to the search of RPOs which, as solutions, have as isotropy subgroup a nontrivial subgroup of $G\times S^1$. Before we get into the statement and proof of this generalization we study in detail the $G\times S^1$--action and its subgroups, and we explain in detail what we mean by nontrivial subgroups. 

All along this section we will assume that the $G$--action on the resonance space $U\sno$ is {\bfi $G$--simple}, that is, $U\sno$ contains a $G$--stable subspace which is either non absolutely irreducible or is isomorphic to the direct sum of two copies of the same absolutely irreducible representation. In the Hamiltonian symmetric framework, this hypothesis occurs generically~\cite[Theorem 3.3]{dellnitz melbourne marsden}. 
Under the $G$--simplicity hypothesis we have the following result
whose proof can be found in~\cite[Proposition 7.2, page 300]{g2}:

\begin{proposition}
\label{spatiotemporal subgroups}
Let $H\subset G\times S^1$ be an isotropy subgroup of the $G\times S^1$--action on the resonance space $U\sno$. Let $\pi:G\times S^1\rightarrow G$ be the projection on the first factor and $K:=\pi(H)\subset G$. If the $G$--action on $U\sno$ is $G$--simple, then:
\begin{enumerate}
\item[{\bf (i)}] $\pi:H\rightarrow K$ is an isomorphism, hence $\dim H=\dim K$.
\item[{\bf (ii)}] There is a homomorphism $\theta\sh:K\rightarrow S^1$ such that 
\[
H=\{(k,\theta\sh(k))\in G\times S^1\mid k\in K\}.
\] 
\item[{\bf (iii)}] $N(H)=N_G(K)\times S^1$.
\end{enumerate} 
\end{proposition}

\noindent Notice that, in the language of the previous proposition, the isotropy subgroups
$K$ of the $G\times S^1$--action on $U\sno$ considered in
Theorem~\ref{theorem 1} are those for which the homomorphism
$\theta\sh$ is identically zero. These are the so called {\bfi
spatial symmetries}. The isotropy subgroups for which $\theta\sh$
is different from zero  are called {\bfi spatiotemporal
symmetries}, and they will be the subject of this section. The
homomorphism $\theta\sh$ will be called the {\bfi temporal
character} of $H$ and its derivative at the identity
$\rho_H:=T_e\theta\sh\in\gk\sus$ the {\bfi temporal velocity} of
$H$. The symbol $\mathfrak{k}$ denotes the Lie algebra of $K$ and
$\mathfrak{k}^\ast $ its dual. The temporal velocity allows us to
express the Lie algebra
$\h$ of $H$ in a very convenient form:
\[
\h=\{(\kappa,T_e\theta\sh(\kappa))\in\gk\times\mathbb{R}\mid\kappa\in\gk\}=\{(\kappa,\langle\rho_H,\kappa\rangle)\in\gk\times\mathbb{R}\mid\kappa
\in\gk\}.
\]

We are now in position to state a generalization of
Theorem~\ref{theorem 1} that incorporates spatiotemporal
symmetries.
 
\begin{theorem}
\label{theorem 1 spatiotemporal}
Let $(V,\omega,h,G,\J:V\rightarrow\gs)$ be a Hamiltonian system with symmetry, with $V$ a vector space, and $G$ a compact positive dimensional Lie group that acts on $V$ in a linear and canonical fashion. Suppose that $h(0)=0$, $\bd h(0)=0$ (that is, the Hamiltonian vector field $X_h$ has an equilibrium at the origin) and that the linear Hamiltonian vector field $A:=DX_h(0)$ is non degenerate and contains $\pm i\no$ in its spectrum. Let $U\sno$ be the resonance space of $A$ with primitive period $T\sno:=\frac{2\pi}{\no}$. Consider the $G\times S^1$--action on $U\sno$, where the $S^1$--action is induced by the semisimple part of $A$, and the Lie group $G$ acts simply on $U\sno$
Let $H=\{(k,\theta\sh(k))\mid k\in K\subset G\}\subset G\times S^1$ be an isotropy subgroup of the $G\times S^1$--action on $U\sno$ with temporal character $\theta\sh$, temporal velocity $\rho\sh\in\gk\sus$, and such that the quadratic form $Q\suh$ on the $H$--fixed point space $U\sno\suh$ defined by
\[
Q\suh(v):=\frac{1}{2}\bd^2 h(0)(v,v),\qquad v\in U\sno\suh
\]
is definite. Then, for any $\chio\in(\gk^\circ)^K$ for which $\J|_{(U\sno)\sh}\inv\left(\chio-\frac{1}{\no}\rho\sh\right)\cap Q_H\inv(1)$
 is non empty ($Q\sh:=Q\suh|_{(U\sno)\sh}$) there exists an open neighborhood $V\schio$ of $\chio$ in $(\gk^\circ)^K$ such that for any $\chi\in V\schio$, the intersection $\J|_{(U\sno)\sh}\inv\left(\chi-\frac{1}{\no}\rho\sh\right)\cap Q_H\inv(1)$ is a submanifold of $(U\sno)\sh$ of dimension $\dim U\sno^H-\dim\ele$.
Suppose that the following  two generic hypotheses hold:
\begin{enumerate}
\item[{\bf (H1)}]  The restriction $h|_{U\sno\suh}$ of the Hamiltonian $h$ to the fixed point subspace $U\sno\suh$ is not radial with respect to the norm associated to $Q\suh$.
\item[{\bf (H2)}] Let  $h_k(v):=\frac{1}{k!}\bd^k h(0)\left(v^{(k)}\right)$, $v\in U\sno\suh$ be the first non radial term in the Taylor expansion of $h|_{U\sno\suh}$ around zero. We will assume that $k\geq 4$ and that the restrictions of $h_k$ to the submanifolds $\J|_{(U\sno)\sh}\inv\left(\chi-\frac{1}{\no}\rho\sh\right)\cap Q_H\inv(1)$, with $\chi\in V\schio$, 
are $\left(N_G(K)_{\rho\sh}\cap N_G(K)_\chi\right)\times
S^1$--Morse.  
\end{enumerate}
Then,  for any  $\epsilon>0$ close enough to zero, $\chi\in V\schio$, and $\lambda:=\Xi\sus(\chi-\frac{1}{\no}\rho\sh,\frac{1}{\no})$, there are at least
\begin{equation}
\label{main estimate 1 spatiotemporal}
{\rm max}\left[\frac{1}{2}\left(\dim U\sno\suh-\dim N_G(K)-\dim \left(N_G(K)_{\rho\sh}\cap N_G(K)_{\chi}\right)+2\dim K\right),\chi_E\left(\EJ\selehd\inv(\lambda)\right)^{L_\lambda}\right]
\end{equation}
distinct relative periodic orbits of $X_h$ with energy $\epsilon$, momentum $\epsilon(\chi-\frac{1}{\no}\rho\sh)\in\gs$, isotropy subgroup $H$, and relative period close to $T\sno$. 
By definition
$\EJ\selehd\inv(\lambda):=
\J|_{(U\sno)\sh}\inv(\chi-s\rho\sh)\cap Q_{A_s}\inv(s\no)
$ and $L_\lambda:=\left(\left(N_G(K)_{s\rho\sh}\cap N_G(K)_{\chi}\right)\times
S^1\right)$.
The symbol $\chi_E\left(\EJ\selehd\inv(\lambda)\right)^{L_\lambda}$
denotes the $L_\lambda$--Euler characteristic of 
$\EJ\selehd\inv(\lambda)$ (which in this case equals the standard
Euler characteristic of the symplectic quotient
$\chi_E(\EJ\selehd\inv(\lambda)/L_\lambda)$).
\end{theorem}

\subsection{Proof of Theorem~\ref{theorem 1 spatiotemporal}} 

In the sequel we will think of the Lie algebra $\gk$ and its dual 
$\gk\sus$ as subspaces of $\g$ and $\g\sus$, respectively, by
choosing in $\g$ an $\Ad_{N_G(K)}$--invariant inner product
$\langle\cdot,\cdot\rangle$ and making $\g=\gk\oplus\m$, with $\m$
the orthogonal complement to $\gk$ in $\g$ with respect to
$\langle\cdot,\cdot\rangle$. If we use the inner product dual to
$\langle\cdot,\cdot\rangle$ we can write $\gs=\gk\sus\oplus\ms$.

We start the proof with the following proposition that provides
several important facts about the temporal velocity and its
relation with the Lie algebras $\mathfrak{h}$ and $\mathfrak{k}$ of
$H$ and $K$, respectively:

\begin{proposition}
\label{properties h}
Let $H\subset G\times S^1$ be an isotropy subgroup of the $G\times S^1$--action on the resonance space $U\sno$, where the $G$--action on $U\sno$ is $G$--simple. Let $K:=\pi(H)\subset G$, $\theta\sh$ be the temporal character of $H$, and $\rho_H:=T_e\theta\sh\in\gk\sus$ be its temporal velocity. Then:
\begin{enumerate}
\item[{\bf (i)}] $\rho\sh\in(\gk\sus)^K$.
\item[{\bf (ii)}] $\h^\circ=\{(-s\rho\sh+\chi,s)\in\gs\times\mathbb{R}\mid s\in\mathbb{R},\,\chi\in\gk^\circ\}=\{(-s\rho\sh+\mu,s)\in\gs\times\mathbb{R}\mid s\in\mathbb{R},\,\mu\in\ms\}$, \\
where $\h^\circ$ denotes the annihilator of $\h$ in $\gs\times\mathbb{R}$ and $\gk^\circ$ that of $\gk$ in $\g\sus$.
\item[{\bf (iii)}] $(\h^\circ)^H=\{(-s\rho\sh+\chi,s)\in\gs\times\mathbb{R}\mid s\in\mathbb{R},\,\chi\in(\gk^\circ)^K\}$
\end{enumerate}
\end{proposition}

\noindent\textbf{Proof \ \ (i)}  It is a consequence of the fact that the temporal character $\theta\sh$ is a homomorphism. Indeed, for any $k\in K$ and $\eta\in\gk$ we have that
\begin{multline*}
\langle\Ad\sus_{k\inv}\rho\sh,\eta\rangle=\langle\rho\sh,\Ad_{k\inv}\eta\rangle
=T_{e}\theta\sh\cdot\Ad_{k\inv}\eta
=\ddto\theta\sh(k\inv\exp t\eta k)=T_e\theta\sh\cdot\eta=\langle\rho\sh,\eta\rangle.
\end{multline*}

\noindent\textbf{(ii)} Recall that we think of $\gk\sus$ as a subspace of $\gs$ by means of the splitting $\gs=\gk\sus\oplus\ms$. By definition:
\begin{align*}
\h^\circ&=\{(\eta,s)\in\gs\times\mathbb{R}\mid\langle\eta,\kappa\rangle+sT_e\theta\sh\cdot\kappa=0, \forall\kappa\in\gk\}\\
	&=\{(\eta,s)\in\gs\times\mathbb{R}\mid\langle\eta+s\rho\sh,\kappa\rangle=0, \forall\kappa\in\gk\}\\
	&=\{(\eta,s)\in\gs\times\mathbb{R}\mid\eta+s\rho\in\gk^\circ\}\\
	&=\{(-s\rho\sh+\chi,s)\in\gs\times\mathbb{R}\mid s\in\mathbb{R}, \chi\in\gk^\circ\}\\
	&=\{(-s\rho\sh+\mu,s)\in\gs\times\mathbb{R}\mid s\in\mathbb{R},\,\mu\in\ms\}.
\end{align*}

\noindent\textbf{(iii)} It is a straightforward consequence of the definition and the use of~\textbf{(i)}. \ \ \ $\blacksquare$

\medskip

We now study the globally Hamiltonian character of the $G\times S^1$--action on the resonance space $U\sno$. The momentum map associated to this action is
\[
\begin{array}{cccc}
\EJ&:U\sno&\longrightarrow&\gs\times{\rm Lie}(S^1)\sus\\
	&v&\longmapsto&\left(\J|_{U\sno}(v),\frac{1}{\no}Q_{A_s}(v)\right),
\end{array}
\]
where $Q_{A_s}$ is the quadratic form $\frac{1}{2}\omega(A_s\cdot,\cdot)$ associated to the semisimple part $A_s$ of  the linearization $A$ of the Hamiltonian vector field $X_h$ at the equilibrium.  As a particular case of what we saw in Section~\ref{Proper actions, fixed 
point sets, slices, and normalizers.} we have that the globally Hamiltonian $G\times S^1$--action on $U\sno$ induces, for each isotropy subgroup $H\subset G\times S^1$ globally Hamiltonian actions of   $L:=\ele$ on $(U\sno)_H$ and $U\sno^H$ with momentum maps $\EJ\selehu:U\sno^H\rightarrow\geles$ and $\EJ\selehd:(U\sno)_H\rightarrow\geles$ associated to these actions given by 
\[
\EJ\selehu(v)=\Xi\sus(\EJ(v)),\qquad\EJ\selehd(v)=\Xi\sus(\EJ(v)),
\]
where $\Xi\sus:(\h^\circ)^H\rightarrow\geles$ is the natural $\ele$--equivariant isomorphism  between the $H$--fixed points in the annihilator of $\h$ in $\gs\times\mathbb{R}$ and the dual of the Lie algebra $\geles$ of $\ele$. 

Note that since the $\ele$--action on $(U\sno)\sh$ is free, the corresponding momentum map $\EJ\selehd$ is a submersion onto its image. Let $\Delta:(\h^\circ)^H\rightarrow (\gk^\circ)^K\times\mathbb{R}$  be the isomorphism defined by $\Delta(-s\rho\sh+\chi,s)=(\chi,s)$.  The mapping $\EJ\snkd:(U\sno)\sh\rightarrow(\gk^\circ)^K\times\mathbb{R}$ defined by $\EJ\snkd:=\Delta\circ\EJ|_{(U\sno)\sh}$ is also a submersion onto its image that more specifically
maps, for any $v\in(U\sno)\sh$, as:
\[
\EJ\snkd(v)=\Delta\circ\EJ(v)=\Delta\left(\J(v),\frac{1}{\no}
Q_{H}(v)\right)=\left(\J(v)+\frac{\rho\sh}{\no}Q_{H}(v),
\frac{1}{\no}Q_{H}(v)\right).
\]
The following proposition justifies the notation utilized in the
statement of the theorem and provides the proof of a few facts that
will be needed later on:

\begin{proposition}
\label{equivalence quotients}
We use the notation introduced in the previous paragraphs. Let 
$\lambda\in\geles$ be an element in the image of $\EJ\selehd$ such
that $\lambda=\Xi\sus(-s\rho\sh+\chi,s)$, for some $s\in\mathbb{R}$
and some $\chi\in(\gk^\circ)^K$. Then:
\begin{description}
\item [(i)] $\EJ\selehd\inv(\lambda)=\EJ\snkd\inv(\chi,s)=\J|_{(U\sno)\sh}\inv(\chi-s\rho\sh)\cap
Q_{H}\inv(s\no)$. $\EJ\selehd\inv(\lambda) $ is a submanifold of
$(U _{\nu_{\circ }})_H $ of dimension
\begin{equation}
\label{dimension level set 1}
\dim U\sno\suh-\dim N_G(K)+\dim K-1.
\end{equation}
\item  [(ii)] The coadjoint
isotropy subgroup $L _\lambda\subset L$
of $\lambda\in
\mathfrak{l}^\ast$ can be
written as $L_\lambda=\left(N_G(K)_{s\rho\sh}\cap
N_G(K)_{\chi}\right)\times S^1/H$, where $N_G(K)_{s\rho\sh}$ and 
$N_G(H)_{\chi}$ are the stabilizers  of
$s\rho\sh\in\gk\sus\subset\gs$ and
$\chi\in(\gk^\circ)^K$  with respect to the coadjoint action of
$N_G(K)$ on $\gs$. 
\item [(iii)]The quotient $\EJ\selehd\inv(\lambda)/L_\lambda$ is a
symplectic manifold of dimension
\[
\dim U\sno\suh-\dim L-\dim L_\lambda=\dim U\sno\suh-\dim N_G(K)-\dim \left(N_G(K)_{s\rho\sh}\cap N_G(K)_{\chi}\right)-2+2\dim K.
\]
\item [(iv)] Let $\chio\in(\gk^\circ)^K$
be an element in $(\gk^\circ)^K$ such that the set $S\schio:=
\J|_{(U\sno)\sh}\inv\left(\chio-\frac{1}{\no}\rho\sh\right)\cap
Q_H\inv(1)$ is non empty and by part {\bf (i)} a submanifold of
dimension~(\ref{dimension level set 1}).  Then, there is a
neighborhood $V\schio$ of $\chio$ in $(\gk^\circ)^K$  such that for
any $\chi\in V\schio$, the set $S_\chi:=
\J|_{(U\sno)\sh}\inv\left(\chi-\frac{1}{\no}\rho\sh\right)\cap
Q_H\inv(1)$ is also nonempty and therefore a manifold of
$(U\sno)\sh$ of dimension~(\ref{dimension level set 1}).
\end{description}
\end{proposition}

\noindent\textbf{Proof} {\bf (i)}  Notice that if
$\lambda=\Xi\sus(-s\rho\sh+\chi,s)$, then
\begin{align}
\label{numerator symplectic}
\J|_{(U\sno)\sh}\inv(\chi-s\rho\sh)\cap Q_{A_s}\inv(s\no)&=\EJ\snkd\inv(\chi,s)=\EJ|_{(U\sno)\sh}\inv(\chi-s\rho\sh,s)=\EJ|_{(U\sno)\sh}\inv((\Xi\sus)\inv(\lambda))\\
	&=\left(\Xi\sus\circ\EJ|_{(U\sno)\sh}\right)\inv(\lambda)=\EJ\selehd\inv(\lambda).\notag
\end{align}
Given that $\EJ\selehd $ is
the  momentum map associated to a free and proper action it is a
submersion and its level sets are submanifolds.
Therefore~(\ref{numerator symplectic})  
has as a corollary that
$\J|_{(U\sno)\sh}\inv(\chi-s\rho\sh)\cap Q_{A_s}\inv(s\no)$ is a
submanifold of $(U\sno)\sh$ of dimension $\dim (U\sno)\sh-\dim\ele$.
Now, as $N(H)=N_G(K)\times S^1$ the expression~(\ref{dimension
level set 1}) follows. 

\noindent {\bf (ii)} Since
$\Xi\sus$ is
$\ele$--equivariant, we have that an arbitrary element
$(l,\theta)H\in\ele$ is actually in the isotropy subgroup
$(\ele)_\lambda$ iff $n\cdot(-s\rho\sh+\chi)=-s\rho\sh+\chi$. Given
that $\rho\sh\in\gk\sus$, $\chi\in\ms$, and $\gk\sus$ and $\ms$ are
$N_G(K)$--stable, then $n\cdot(-s\rho\sh+\chi)=-s\rho\sh+\chi$ iff
$n\cdot(-s\rho\sh)=-s\rho\sh$ and $n\cdot\chi=\chi$, that is, iff
$n\in N_G(K)_{s\rho\sh}\cap N_G(K)_{\chi}$. Therefore,
\[
\left(\ele\right)_\lambda=\left(N_G(K)_{s\rho\sh}\cap N_G(K)_{\chi}\right)\times S^1/H.
\]

\noindent {\bf (iii)} It is a straightforward consequence of the
Marsden--Weinstein Reduction Theorem~\cite{mwr}  and points {\bf
(i)} and {\bf (ii)}.

\noindent {\bf (iv)} Let
$\lo=\Xi\sus\left(\chio-\frac{1}{\no}\rho\sh,\frac{1}{\no}\rho
\sh\right)\in\geles$. 
By {\bf (i)},
$S\schio=\EJ\inv\selehd(\lo)$. Since $\EJ\selehd$ is the momentum
map corresponding to a free action on $(U\sno)\sh$ and therefore a
submersion onto its image,  
the mapping $T_v\EJ\selehd:T_v(U\sno)\sh\rightarrow\geles$ is
surjective, for any $v\in S\schio=\EJ\inv\selehd(\lo)$. The
Local Onto Theorem (see for instance~\cite[Theorem 3.5.2]{mta})
implies the existence of an open neighborhood $W\slo$ of $\lo$ in
$\geles$ and an open neighborhood $U_v$ of $v$ in
$(U\sno)\sh$ such that the mapping
$\EJ\selehd|_{U_v}:U_v\rightarrow W\slo$ is onto. In particular,
for any $\lambda\in W\slo$, the level set $\EJ\selehd\inv(\lambda)$
is non empty and, by the submersion argument is a submanifold of
$(U\sno)\sh$ of dimension~(\ref{dimension level set 1}). 
We now construct the open neighborhood $V\slo$ of $\lo$ whose
existence we claim in the statement of the Lemma. Firstly, the set
$T$ defined by
$T:=\Delta\circ(\Xi\sus)\inv(W\slo)\subset(\gk^\circ)^K\times\mathbb{R}$
is an open neighborhood of $(\chio,\frac{1}{\no}\rho\sh)$. By the
openness of $T$, there exist open neighborhoods $V\schio\subset
(\gk^\circ)^K$ and $W_{\frac{1}{\no}\rho\sh}\subset\mathbb{R}$, of
$\chio$ and $\frac{1}{\no}\rho\sh$, respectively, such that
$V\schio\times W_{\frac{1}{\no}\rho\sh}\subset T$. $V\schio$ is the
neighborhood needed in the statement.\ \ \ \ $\blacksquare$

\medskip

We now prove a result that constitutes the spatiotemporal analog of Lemma~\ref{smooth branches}. Let $\widehat{h|_{U\sno\suh}}$ be the equivalent Hamiltonian in normal form associated to $h|_{U\sno\suh}$.

\begin{lemma}
\label{smooth branches spatiotemporal}
Suppose that we are under the hypotheses of Theorem~\ref{theorem 1 spatiotemporal}. Then, the restriction of the function $h_k\in C\suinf(U\sno\suh)$ defined by $h_k(u):=\frac{1}{k!}\bd^k h|_{U\sno\suh}(0)(u^{(k)})$,  to the level sets of the form $\J|_{(U\sno)\sh}\inv\left(\chi-\frac{1}{\no}\rho\sh\right)\cap Q_H\inv(1)$, where $\chi$ sits in $V\schio$, the neighborhood of $\chio\in(\gk^\circ)^K$ introduced in the previous Lemma, has at least
\begin{equation}
\label{critical lusternik spatiotemporal}
{\rm max}\left[\frac{1}{2}\left(\dim U\sno\suh-\dim N_G(K)-\dim \left(N_G(K)_{\rho\sh}\cap N_G(K)_{\chi}\right)+2\dim K\right),\chi_E\left(\EJ\selehd\inv(\lambda)\right)^{L_\lambda}\right]
\end{equation}
distinct critical orbits, where $\lambda:=\Xi\sus(\chi-\frac{1}{\no}\rho\sh,\frac{1}{\no})$ and $\chi_E\left(\EJ\selehd\inv(\lambda)\right)^{L_\lambda}$ denotes the $L_\lambda$--Euler characteristic of  $\EJ\selehd\inv(\lambda)$. 

Furthermore, let $\chi'\in V\schio\subset(\gk^\circ)^K$ be arbitrary but fixed and let $\{u_1,\ldots,u_k\}$ be the set of critical points of the restriction of the function $h_k$ to the level set $\J|_{(U\sno)\sh}\inv\left(\chi'-\frac{1}{\no}\rho\sh\right)\cap Q_H\inv(1)$, provided by~(\ref{critical lusternik spatiotemporal}). Then, for each $u_i$, $i\in\{1,\ldots,k\}$, there exist a neighborhood $E\subset\mathbb{R}$ of the origin in the real line, a neighborhood $V_{\chi'}\subset V\slo\subset (\gk^\circ)^K$ of $\chi'$ in $(\gk^\circ)^K$, and a smooth function $\rho:E\times V_{\chi'}\rightarrow Q\sh\inv(1)$ such that $\rho(0,\chi')=u_i$ and also, the function $v:E\times V_{\chi'}\rightarrow U\sno\suh$ defined by $v(r,\chi):=r\rho(r,\chi)$ satisfies that:
\begin{enumerate}
\item[{\bf (i)}] $v(r,\chi)\in (U\sno)\sh$ iff $r\neq 0$.
\item[{\bf (ii)}] $Q\suh(v(r,\chi))=r^2$ and $\J(v(r,\chi))=r^2(\chi-\frac{1}{\no}\rho\sh)$, $r\in E$, $\chi\in V_{\chi'}$.
\item[{\bf (iii)}] $\bd \widehat{h|_{U\sno\suh}}|_{\J|_{(U\sno)\sh}\inv\left(r^2(\chi-\frac{1}{\no}\rho\sh)\right)\cap Q_H\inv(r^2)}(v(r,\chi))=\bd \widehat{h|_{U\sno\suh}}|_{\EJ\snkd\inv(r^2\chi,\frac{r^2}{\no})}(v(r,\chi))=0$, that is, the branches $v(r,\chi)$ are made of critical points of the restriction of the function $\widehat{h|_{U\sno\suh}}$ to the 
level sets
\[
\J|_{(U\sno)\sh}\inv\left(r^2(\chi-\frac{1}{\no}\rho\sh)\right)\cap
Q_H\inv(r^2).
\]
\end{enumerate}
\end{lemma}

\noindent\textbf{Proof} The estimate~(\ref{critical lusternik spatiotemporal}) is a straightforward consequence of  Proposition~\ref{equivalence quotients}, Corollary~\ref{cup length symplectic}, the statement~(\ref{euler estimate}), and the invariance properties of $h_k$.

We now construct the branches $v(r,\chi)$ in a fashion similar to Lemma~\ref{smooth branches}. Let $\chi'\in V\schio\subset(\gk^\circ)^K$ be arbitrary but fixed and let $u\szero$ be one  of the critical points of the restriction of the function $h_k$ to the level set $\J|_{(U\sno)\sh}\inv\left(\chi'-\frac{1}{\no}\rho\sh\right)\cap Q_H\inv(1)$, provided by~(\ref{critical lusternik spatiotemporal}). Since by 
Proposition~\ref{equivalence quotients} any element of the form
$(\chi,\frac{1}{\no})$ with $\chi\in V\schio$ is a regular value of
$\EJ\snkd$, there exist neighborhoods $V_{\chi'}\subset V\schio$ of
$\chi'$, $W_{\frac{1}{\no}}$ of $\frac{1}{\no}$ in $\mathbb{R}$, 
$W$  of the origin in $\mathbb{R}^s$, with $s=\dim
U\sno^H-\dim(\gk^\circ)^K-1$, and a mapping $\psi:V_{\chi'}\times
W\times W_{\frac{1}{\no}}\rightarrow U\sno^H$ which is a
diffeomorphism onto its image. The mapping $\psi$ satisfies that
$\psi(\chi',0,\frac{1}{\no})=u\szero$ and $\EJ\snkd(\psi(\chi,w,l))=(\chi,l)$, which is equivalent to having that $\J(\psi(\chi,w,l))=\chi-l\rho\sh$ and $Q\sh(\psi(\chi,w,l))=l\no$. Let $\varphi:V_{\chi'}\times W\rightarrow Q\sh\inv(1)$ be the mapping defined by $(\chi,w)\longmapsto\psi(\chi,w,\frac{1}{\no})$. This map constitutes a local chart around the point $u\szero$ in $Q\sh\inv(1)$. Exactly as  we did in Lemma~\ref{smooth branches}, we can factor out the $(N_G(K)_{\rho\sh}\cap N_G(K)_{\chi'})\times S^1$--action in this chart which allows us to prove the Lemma by mimicking what we did in the proof of  Lemma~\ref{smooth branches} starting from expression~(\ref{j and g}).\ \ \
$\blacksquare$

\medskip

Let now $\chi'\in V\schio$ and $v:E\times V_{\chi'}
\rightarrow U\sno\suh$ be one of the branches 
introduced in Lemma~\ref{smooth branches spatiotemporal}.
By construction, 
\[
\widehat{h|_{U\sno\suh}}|_{\J|_{(U\sno)\sh}\inv\left(r^2(\chi-\frac{1}{\no}\rho\sh)\right)\cap Q_H\inv(r^2)}(v(r,\chi))=\bd \widehat{h|_{U\sno\suh}}|_{\EJ\snkd\inv(r^2\chi,\frac{r^2}{\no})}(v(r,\chi))=0.
\] 
Since $\EJ\snkd$ maps into $(\gk^\circ)^K\times \mathbb{R}$, the Lagrange Multipliers Theorem guarantees the existence of an element $\Lambda(r,\chi)\in \left((\gk^\circ)^K\right)\sus$ and $c(r,\chi)\in\mathbb{R}$ such that
\begin{equation}
\label{almost pairing}
\bd \widehat{h|_{U\sno\suh}}(v(r,\lambda))=\left[c(r,\lambda)+\langle\rho\sh,\Lambda(r,\lambda)\rangle\right]\bd Q^H(v(r,\lambda))+\bd(\J|_{U\sno\suh})^{\Lambda(r,\lambda)}(v(r,\lambda)),
\end{equation}
which implies that $v(r,\lambda)$ is a periodic point of the Hamiltonian vector field $X_{\widehat{h|_{U\sno^H}}-(\J|_{U\sno\suh})^{\Lambda(r,\lambda)}}$. If we are able to prove that $\Lambda(r,\lambda)$ becomes very small as $r\rightarrow 0$, the Normal Form Theorem will guarantee that $v(r,\lambda)$ will amount to a periodic point of  $X_{h|_{U\sno\suh}-(\J|_{U\sno\suh})^{\Lambda(r,\lambda)}}=X_{h|_{U\sno\suh}}-(\Lambda(r,\lambda))_{U\sno\suh}$, that is, a RPO of $X_h$ (the symbol $(\Lambda(r,\lambda))_{U\sno\suh}$ denotes the vector field defined by $(\Lambda(r,\lambda))_{U\sno\suh}(v)=\ddto\exp t\Lambda(r,\lambda)\cdot v,\,v\in U\sno\suh$. It is easy to show that this is a well defined vector field on $U\sno\suh$ since  $\exp t\Lambda(r,\lambda)\cdot v\in U\sno\suh$ whenever $v\in U\sno\suh$). Actually it can be easily proved by mimicking what we did in the proof of Theorem~\ref{theorem 1} after~(\ref{identity pairing}), that $\Lambda(r,\chi)$ and $c(r,\chi)$ are smooth functions that tend to zero and one, respectively, as the variable $r$ tends to zero.\ \ \ $\blacksquare$

\medskip

\begin{remark}
\label{difference between theorems}
\normalfont
As we already pointed out in the introduction Theorem~\ref{theorem
1 spatiotemporal} is NOT a generalization of Theorem~\ref{theorem
1}. In the proof of Theorem~\ref{theorem 1} intervenes a
transversality argument that guarantees that all the solutions
obtained are genuine RPOs and not just relative equilibria. This is
explicitly mentioned in its statement.
Even though the spatial symmetries in the relative periodic
solutions predicted in Theorem~\ref{theorem 1} are a particular
case of the  spatiotemporal symmetries treated in 
Theorem~\ref{theorem 1 spatiotemporal} one does not generalize the
other since the subgroups of
$G \times S ^1 $ intertwine the $G$  and $S ^1 $--actions via the
temporal character preventing us from making the distinction
between RPOs and relative equilibria. In short, we cannot guarantee
that the RPOs provided by Theorem~\ref{theorem 1
spatiotemporal} are not just relative equilibria.
\ \ \ $\blacklozenge$
\end{remark}

\section{Relative periodic orbits around stable relative equilibria}
\label{relative periodic orbits around stable relative equilibria}

In this section we will use the so called
Marle--Guillemin--Sternberg (MGS) normal form and the 
reconstruction equations in order to generalize the main result in
the previous section to the search of RPOs around genuine relative
equilibria.

\subsection{The MGS normal form and the reconstruction equations}

Since this topic has been already introduced already in many other 
papers  we will just briefly sketch the results that we will need
in our exposition, and will leave the reader interested in the
details consult the original papers~\cite{nfm, nfgs}. Regarding the
reconstruction equations the reader is encouraged to check
with~\cite{thesis, roberts wulff 99, review}.

All along this section we will work with a $G$--Hamiltonian system  $(M,\omega,h,G,\J)$, where the Lie group $G$ acts in a proper and globally Hamiltonian fashion on the manifold $M$. Let $m$ be a point in $M$ such that $\J(m)=\mu\in\gs$ and
$G_m$ denotes the isotropy subgroup of the point $m$. We denote
by $\gmu$ the Lie algebra of the stabilizer $G\sm$ of $\mu\in\gs$
under the coadjoint action of $G$ on
$\gs$.  We now choose in 
$\ker T\seme \J$ a $G_{m}$--invariant inner product, $\langle\cdot,\cdot\rangle$, always available by the
compactness of $G_{m}$. Using this inner product we define the {\bfi symplectic normal space}  $V_m$ at $m\in M$ with respect to the inner product $\langle\cdot,\cdot\rangle$, as the
orthogonal complement of $T\seme(G\sm\cdot m)$ in $\ker T\seme\J$, 
that is,
$
\ker T\seme\J=
T\seme(G\sm\cdot m)\oplus V_m$,
where the symbol $\oplus$ denotes orthogonal direct sum.
It is easy to verify that $(V_m,\,\omega(m)|_{V_m})$ is a $G_{m}$--invariant symplectic
vector space. Let $B ^{\sharp}: V _m^\ast  \rightarrow V _m $ be the
isomorphism associated to the symplectic form
$\omega(m)|_{V_m} $

Recall that by the equivariance of $\J$, the isotropy
subgroup $G_{m}$ of $m$ is a subgroup of $G\sm$ and therefore
$\mathfrak{g}_{m}={\rm Lie}(G_{m})\subset\gmu$. Using again the compactness of $G_{m}$, 
we construct an inner product
$\langle \cdot,\cdot \rangle $ on $\g$, invariant under the
restriction to $G_{m}$ of the adjoint action of $G$ on $\g$. Relative
to this inner product we can write  the following orthogonal 
direct sum decompositions
$
\mathfrak{g}=\mathfrak{g}_{\mu}\oplus
\mathfrak{q}$, and
$\mathfrak{g}_{\mu}=\mathfrak{g}_{m}\oplus \mathfrak{m}$,
for some subspaces $\mathfrak q \subset \g$ and $\mathfrak m \subset
\mathfrak g_\mu$. The inner product also allows us to identify all 
these Lie algebras with their duals. In particular, we have the 
dual orthogonal  direct sums $\gs=\gmus\oplus\mathfrak q \sus$
and $\gmus=\mathfrak{g}_m^\ast\oplus\ms$ which allow us to consider 
\gmus\ as a subspace of \gs\, and, similarly, $\mathfrak{g}_m^\ast$
and
\ms\ as subspaces of $\gmus$. 

The $G_{m}$--invariance of the inner product utilized to construct 
the splittings $\mathfrak{g}_{\mu}=\mathfrak{g}_{m}\oplus \mathfrak{m}$ and
$\gmus=\mathfrak{g}_m^\ast\oplus\ms$, implies that  both $\m$ and
$\ms$ are
$G_{m}$--spaces using the restriction to them of the $G_{m}$--adjoint and
coadjoint actions, respectively. 

The importance of all these objects is in the fact that 
there
is a positive number $r>0$ such that, denoting by $\ms_r$ the open
ball of radius $r$ relative to the $G_{m}$--invariant inner product on
$\ms$, the manifold
$Y_r:=G\times_{G _m}(\ms_r\times V_m)$
can be endowed with a symplectic structure $\omega_{Y_r}$ with 
respect to which the left $G$--action
$g\cdot[h,\,\eta,\,v]=[gh,\,\eta,\,v]$ on $Y_r$ is globally
Hamiltonian with ${\rm Ad}^*$--equivariant momentum map
$\J_{Y_r}:Y_r \rightarrow \gs$  given by
$
\J_{Y_r}([g,\rho,v])={\rm Ad}\sus_{g\inv}\cdot
(\mu+\rho+\J_{V_m}(v))$.
Moreover, 
there exist $G$--invariant
neighborhoods $U$ of $m$ in $M$, $U'$ of $[e,\,0,\,0]$ in $Y_r$, and
an equivariant symplectomorphism $\phi:U\rightarrow U'$ satisfying
$\phi(m)=[e,\,0,\,0]$ and $\J_{Y_r}\circ\phi=\J$. On other words,
the twisted product $Y _r$ can be used as a coordinate system in a
tubular neighborhood of the orbit $G \cdot m $. This semi--global
coordinates are referred to as the MGS normal form.

In what follows we
will use the MGS coordinates to compute the equations that describe
the dynamics induced by the  Hamiltonian vector field corresponding
to a
$G$--invariant Hamiltonian. These are called the
bundle~\cite{roberts wulff 99} or  reconstruction~\cite{thesis} 
equations. Let
$h\in C\suinf(Y)^G$ be a
$G$--invariant Hamiltonian on $Y$.  Our aim is to compute the
differential equations that determine the $G$--equivariant Hamiltonian vector field
$X_h\in \mathfrak{X}(Y)$ associated to $h$ and characterized by 
$
\mathbf{i}_{X_h}\omega_Y=\bd h$.

Since the projection $\pi:G\times\ms\times V_m\rightarrow 
G\times\sh(\ms\times V_m)$ is a surjective submersion, there are
always local sections available  that we can use to locally express
$X_h=T\pi (X_G,\,X_{\ms},\,X_{V_m})$,
with $X_G,\,X_{\ms}$ and $X_{V_m}$ locally defined smooth maps
on $Y$ and having values in $TG,\,T\ms$ and $TV_m$
respectively. Thus, for any $[g,\,\rho,\,v]\in Y$, one has 
$X_G([g,\,\rho,\,v])\in T_gG$, $X_{\ms}([g,\,\rho,\,v])\in
T_\rho\ms = \ms$, and $X_{V_m}([g,\,\rho,\,v])\in T_v{V_m} = V_m$.
Moreover, using the $\Ad_{G _m}$--invariant decomposition of the Lie
algebra $\g$:
$\g=\mathfrak{g}_{m}\oplus\m\oplus\q$,
the mapping $X_G$ can be 
written, for any $[g,\,\rho,\,v]\in Y$, as
$
X_G([g,\,\rho,\,v])=T_eL_g\big(X_{\mathfrak{g}_{m}}([g,\,\rho,\,v])+X_{\m}([g,\,
\rho,\,v])+X_{\q}([g,\,\rho,\,v])\big)$,
with $X_{\mathfrak{g}_{m}}$, $X_{\m}$, and $X_{\q}$, locally
defined smooth maps on $Y$ with values in $\mathfrak{g}_{m},\,\m$,
and $\q$ respectively. Also, note that since $h\in C\suinf(G\times\sh(\ms\times V_m))^G$ is 
$G$--invariant, the mapping $h\circ\pi\in C\suinf(G\times\ms\times
V_m)^H$ can be understood as a $H$--invariant function that depends
only on the $\ms$ and $V_m$ variables, that is,
$
h\circ\pi\in C\suinf(\ms\times V_m)\suh $.
 
Using these ideas and the explicit expression of the symplectic 
form $\omega_{Y_r}$ we can
explicitly write down the differential equations that determine the
components of $X_h$. In order to do so we first implicitly define a
function $\eta:\gmu\times\q\rightarrow\qs$ such that
$\eta(\xi,0)=0$ for all $\xi\in\gmu$ and that for $\rho\in\ms$,
$v\in V_m$ small enough satisfies
$
\mathbb{P}_{\q\sus}\left(\ad\sus_{(X_{\m}+\eta(X_{\m},\rho+\J_{V_m}(v)))}
      (\rho+\J_{V_m}(v)+\mu)\right)=0$.
If we define
$
\psi(\rho,v):=\eta(D_{\ms}(h\circ\pi),\rho+\J_{V_m}(v))$, then:
\medskip
\begin{center}
\begin{boxit}
\begin{gather}
X_{\mathfrak{g}_{m}}=0\\
X_\q=\psi(\rho,v)\label{field 0}\\
X_{\m}=D_{\ms}(h\circ\pi)
\label{field 1 general}\\
X_{V_m}=B^\sharp_{V_m}(D_{V_m}(h\circ\pi))\label{field 2 general}\\
X_{\ms}=\mathbb{P}_{\m\sus}\left(\ad\sus_{D_{\ms}(h\circ\pi)}
\rho+\ad\sus_{D_{\ms}(h\circ\pi)}\J_{V_m}(v)+
\ad\sus_{\psi(\rho,v)}(\rho+\J_{V_m}(v))\right).\label{field
3 general}
\end{gather}
\end{boxit}
\end{center}

\medskip

The previous equations admit severe simplifications in the presence of  various Lie algebraic hypotheses. See~\cite{roberts wulff 99} for an extensive study. For future reference we will note two particularly important cases:
\begin{itemize}
\item The Lie algebra $\g$ is Abelian: in that case
$X_{\mathfrak{m}^\ast}=X_\q=0$ at any point.
\item The point $\mu\in\gs$ is {\bfi split}~\cite{lerman book}, that is, the Lie algebra $\gmu$ of the coadjoint isotropy of $\mu$ admits a ${\rm Ad}_{G\sm}$--invariant 
complement in $\mathfrak{g}$: in that case the mappings $\eta$ and
$\psi$ are identically zero.
\end{itemize}

\subsection{The main estimate}

The following result generalizes theorems~\ref{theorem 1} and~\ref{theorem 1 spatiotemporal} to the search of  RPOs around stable relative equilibria

\begin{theorem}
\label{rpos around stable re}
Let $(M,\omega,h,G,\J)$ be a Hamiltonian $G$--space. Let $m\in M$ be a relative equilibrium of this system with velocity $\xi\in\g$,
isotropy $G\seme$, $\J(m)=\mu\in\gs$, and $h(m)=0$. Let
$V\seme\subset T\seme M$ be any symplectic normal space through the
point $m$,
$Q_{V\seme}:=\frac{1}{2}\bd^2\left(h-\J^{\mathbb{P}_\m\xi}\right)(m)|_{V\seme}$,
and $A_{V\seme}:=X_{Q_{V\seme}}$. Suppose that for $V\seme$  (and
hence for any other symplectic normal space) the infinitesimally
symplectic linear map $A_{V\seme}$ is non singular and has $\pm
i\no$ as eigenvalues. Let $U\sno$ be the $G_{m}\times S^1$--symmetric
resonance space of $A_{V\seme}$ with primitive period $T\sno$ (the
$S^1$--action is generated by the semisimple part of the
restriction of $A_{V\seme}$ to $U\sno$).  Let $H\subset G_{m}\times
S^1$ be an isotropy subgroup of the $G_{m}\times S^1$--action on
$U\sno$. If $H$ is not a purely spatial subgroup of $G_{m}\times S^1$
we assume that the $G_{m}$--action on $U\sno$ is simple, and hence it
has an associated temporal character $\theta_{H}$ such that
\[
H=\{(f,\theta_{H}(f))\mid f\in K\subset G_{m}\}
\]
and a well defined temporal velocity
$\rho_{H}:=T_e\theta_{H}\in\frak{k}\sus$.

Suppose that
$\bd^2\left(h-\J^{\mathbb{P}_\m\xi}\right)(m)|_{U\sno^{H}}$ is a
definite quadratic form. We consider two cases:
\begin{enumerate}
\item[{\bf (i)}] If $\rho_{H}=0$ and
$\J_{V\seme}|_{(U\sno)_{H}}\inv(0)\cap Q_{V\seme}\inv(1)$ is
nonempty then,  for any $\epsilon>0$ small enough, there exist
generically
\begin{equation}
\label{main estimate 1 relative 1}
{\rm max}\left[\frac{1}{2}\left(\dim U\sno^{H}-2\dim
N_{G_{m}}(K)+2\dim
K\right),\chi_E\left(\J_{V\seme}|_{(U\sno)_{H}}\inv(0)\cap
Q_{V\seme}\inv(1)\right)^{N_{G_{m}}(K)\times S^1}\right]
\end{equation}
distinct RPOs of $X_h$ with energy $\epsilon$, momentum $\mu\in\gs$, relative period close to $T\sno$, and isotropy subgroup $G_{m}$.
\item[{\bf (ii)}] If $\rho_{H}\neq 0$ and
$\chio\in(\frak{k}^\circ)^K$ is such that 
$\J_{V\seme}|_{(U\sno)_{H}}\inv\left(\chio-\frac{1}{\no}\rho\sh\right)\cap
Q_{V\seme}\inv(1)$
 is non empty, assume that ONE of the following hypotheses holds:
\begin{enumerate}
\item[{\bf 1.}] The Lie algebra $\g$ is Abelian.
\item[{\bf 2.}] The Lie algebra $\gmu$ is Abelian and $\mu$ is split.
\item[{\bf 3.}] $\mathfrak{g}_{m}=\gmu$.
\end{enumerate}
Then, there exists an open neighborhood $V\schio$ of $\chio$ in
$(\frak{k}^\circ)^K$ such that for any $\chi\in V\schio$, the
intersection
$\J_{V\seme}|_{(U\sno)_{H}}\inv\left(\chi-\frac{1}{\no}\rho\sh\right)\cap
Q_{V\seme}\inv(1)$ is a submanifold of $(U\sno)_{H}$ of dimension
$\dim U\sno^{H}-\dim N_{G_{m}\times S^1}(H)/H$. Moreover, for any 
$\epsilon>0$ close enough to zero and $\chi\in V\schio$, there are
generically at least
\begin{multline}
\label{main estimate 1 spatiotemporal relative}
{\rm max}\left[\frac{1}{2}\left(\dim U\sno^{H}-\dim
N_{G_{m}}(K)-\dim \left(N_{G_{m}}(K)_{\rho_{H}}\cap
N_{G_{m}}(K)_{\chi}\right)+2\dim K\right)\right.,\\
\left.\chi_E\left(\J_{V\seme}|_{(U\sno)_{H}}\inv(\chi-\frac{1}{\no}\rho_{H})\cap
Q_{V\seme}\inv(1)\right)^{\left(N_{G_{m}}(K)_{\rho\sh}\cap
N_{G_{m}}(K)_{\chi}\right)\times S^1}\right]
\end{multline}
distinct relative periodic orbits of $X_h$ with energy $\epsilon$,
momentum 
$\mu+\epsilon(\chi-\frac{1}{\no}\rho_{H})\in\gs$,  isotropy $H$, and
relative period close to $T\sno$.
\end{enumerate}

The symbol $\J_{V\seme}$ denotes the momentum map associated to the linear action of $G_{m}$ on $V\seme$. The projections
$\mathbb{P}_{\mathfrak{g}_{m}}$ and $\mathbb{P}_\m$ are consistent with a given
${\rm Ad}_{G_{m}}$--invariant splitting $\g=\mathfrak{g}_{m}\oplus\m\oplus\q$ of
the  Lie algebra $\g$.
\end{theorem}

\begin{remark}
\normalfont
The reader interested in the relation between the hypothesis on the definiteness of the quadratic form $\bd^2\left(h-\J^{\mathbb{P}_\m\xi}\right)(m)|_{V\seme}$ and the actual nonlinear stability of the relative equilibrium is encouraged to check with~\cite{singreleq,
re, patrick roberts wulff 2002}, and references therein.\ \ \
$\blacklozenge$
\end{remark}

\noindent\textbf{Proof of the Theorem} A straightforward
computation shows that  the Hessian in the statement of the
theorem is well defined and that the hypotheses on it do not depend
on the choice of symplectic normal space $V_m$. 

Given the local nature of the statement, we can use the MGS
coordinates to carry out the proof of the theorem. For simplicity
in the exposition we will identify points and maps in $M$ and their
counterparts in the MGS coordinates $Y$. Those coordinates can be
chosen so that the point $m$ is represented by $[e,0,0]\in
G\times_{G_{m}}(\ms\times V\seme)$ and the submanifold
$\Sigma\seme:=\{e\}\times_{G_{m}}(\{0\}\times V\seme)\subset Y$ is
such that $T\seme\Sigma\seme$ is a symplectic normal space at $m$,
that is, $\ker T\seme\J=T\seme\Sigma\seme\oplus T\seme(G\sm\cdot
m)$. 

Notice that in MGS coordinates the point $m\equiv[e,0,0]$ is a relative equilibrium of the Hamiltonian vector field $X_h$ with velocity $\xi$ when 
\begin{equation}
\label{equilibrium}
X_h(m)=T_{(e,0,0)}\pi (\xi,0,0)
\end{equation}
and the associated flow behaves as $F_t(m)=[\exp t\xi,0,0]$.

We now define the function $h_{V\seme}\in C\suinf(V\seme)^{G_{m}}$ as
$h_{V\seme}(v)=h\circ\pi(0,v)$, for each $v\in V\seme$. Moreover,
notice that by~(\ref{equilibrium}) and the reconstruction
equation~(\ref{field 2 general}) 
\[
\bd h_{V\seme}(0)=D_{V\seme}(h\circ\pi)(0,0)=B^\flat_{V\seme}(X_{V\seme}(0,0,0))=0,
\]
where $B\in\Lambda^2(\vm\times\vm)$ is the Poisson tensor associated to the symplectic form $\omega\svm:=\omega|\svm$.
Also, for any $v, w\in V\seme$:
\begin{multline*}
\bd^2(h-\J^{\mathbb{P}_\m\xi})([e,0,0])(T_{(e,0,0)}\pi(0,0,v),T_{(e,0,0)}\pi(0,0,w))=\ddto\ddso(h-\J^{\mathbb{P}_\m\xi})([e,0,tv+sw])\\
=\bd^2 h_{V\seme}(0)(v,w)-\ddto\langle T_{tv} \J_{V\seme}\cdot w,\mathbb{P}_\m\xi\rangle=\bd^2 h_{V\seme}(0)(v,w),
\end{multline*}
since $T_{tv} \J_{V\seme}\cdot w\in\mathfrak{g}_{m}^\ast$ for any
$t$. Therefore, the hypothesis on the non degeneracy of
$\bd^2(h-\J^{\mathbb{P}_\m\xi})(m)|_{U\sno}$ implies the non
degeneracy of $\bd^2 h_{V\seme}(0)|_{U\sno}$.

We now apply Theorem~\ref{theorem 1 spatiotemporal} to the equilibrium that the system $(V\seme,\omega_{V\seme},h_{V\seme},G_{m},\J_{V\seme})$ has at the origin. If we use isotropy subgroups of the $G_{m}\times S^1$--action on $V\seme$ with temporal velocity equal to zero and look for RPOs such that $\chi=0$, Theorem~\ref{theorem 1 spatiotemporal} provides us with~(\ref{main estimate 1 relative 1}) RPOs for $(V\seme,\omega_{V\seme},h_{V\seme},G_{m},\J_{V\seme})$ with $\J_{V\seme}$ momentum equal to zero. The general case with non zero temporal character and arbitrary $\chi$ gives us~(\ref{main estimate 1 spatiotemporal relative}) relatively periodic solutions with $\J_{V\seme}$ momentum equal to
$\epsilon(\chi-\frac{1}{\no}\rho_{H})$.

In the remainder of the proof we will use these $G_{m}$--relative periodic orbits in $\vm$ to construct $G$--relative periodic orbits in the original system using the reconstruction equations. We will first establish the estimate~(\ref{main estimate 1 relative 1}) on the number of RPOs with momentum equal to $\mu$: let $v\in\vm$ be one of  the  $G_{m}$--relative periodic orbits of $(V\seme,\omega_{V\seme},h_{V\seme},G_{m},\J_{V\seme})$ with $\J_{V\seme}$ momentum equal to zero. If we look at the reconstruction equation~(\ref{field 3 general}) taking into account that $\J_{\vm}(v)=0$ we obtain that the point $[e,0,v]$ is such that $X_\ms([e,0,v])=0$ and therefore it is necessarily a $G$--relative periodic point of $X_h$. Notice that, given the    expression of the momentum map in MGS
coordinates, this RPO has momentum exactly equal to
$\mu$. 

As to the estimate~(\ref{main estimate 1 spatiotemporal relative}), consider now one of the RPOs of $(V\seme,\omega_{V\seme},h_{V\seme},G_{m},\J_{V\seme})$ with 
$\J_{V\seme}$ momentum equal to
$\epsilon(\chi-\frac{1}{\no}\rho_{H})$. Additionally, suppose that
we are in  any of the first two cases contemplated in the Lie
algebraic hypotheses in the statement of the theorem, that is, 
either the Lie algebra $\g$ is Abelian or $\gmu$ is Abelian and
$\mu$ is split. It is easy to see by looking at the reconstruction
equation~(\ref{field 3 general}) that in any of those two cases
$X_\ms=0$ at any point and therefore if $v\in V\seme$ is one of the
$G_{m}$--RPOs of
$(V\seme,\omega_{V\seme},h_{V\seme},G_{m},\J_{V\seme})$  the point
$[e,0,v]$ is necessarily a $G$--RPO of the original system, with
$G$--momentum map $\mu+\epsilon(\chi-\frac{1}{\no}\rho_{H})\in\gs$
and isotropy subgroup $H$. If we are under hypothesis 3, the fact
that $\mathfrak{g}_{m}=\gmu$ implies that $\ms=0$ and therefore the
argument that we just used can be applied to the points of the form
$[e,v]$.
\ \ \
$\blacksquare$

\medskip

\addcontentsline{toc}{section}{Conclusions}
\noindent\textbf{Conclusions. } In this paper we have proved
results that give estimates on the number of relative periodic
orbits around given stable equilibria and relative equilibria.

The approach taken in the
proofs carries in its wake some limitations in our results.  For
instance, discrete symmetries are invisible by the momentum map. It
is our belief that results in this direction can only be obtained
by taking a global variational approach that the author is already
studying and that will be the subject of a future work.

This global variational approach seems also to be the only way to obtain global generalizations of the local results stated throughout the paper, similar to those obtained in the past regarding the Weinstein--Moser Theorem (see for instance~\cite{ekeland lasry 80, ambrosetti mancini 82, hofer zehnder 95}, an references therein)  where, by substituting the stability condition  by convexity hypotheses, estimates regarding the existence of periodic orbits could be formulated for any convex energy level set.

\bigskip

\addcontentsline{toc}{section}{Acknowledgments}
\noindent\textbf{Acknowledgments.}  I am thankful to J. Marsden for emphasizing to 
me the importance of the so called \emph{Smale Program} that, at some level, played
an important inspirational role in this work. Special thanks go also to J. Montaldi
for his help and for his organizational skills that made the summer school MASESS
2000 at Peyresq (France) a comfortable place where  I developed some ideas that
improved this paper. I also thank A. Weinstein for kindly answering my questions,
B. \^Zhilinski\'{\i} for his interest and suggestions, and T. Ratiu for his
encouragement. The valuable comments of an anonymous referee are
also gratefully acknowledged since they had a very positive impact
in this paper. This
research was partially supported by the European Commission through
funding for the Research Training Network
\emph{Mechanics and Symmetry in Europe} (MASIE).

\end{document}